\magnification=1200
\magnification=1200
\hsize=6 true in
\vsize=7.5 true in
\voffset=.5 true in
\tolerance=1600

\font\rmtwelve=cmbx10 at 12pt

\vskip 3pt
\bigbreak
\bigbreak
\bigbreak

\centerline{\rmtwelve Family Seiberg-Witten invariants 
and wall
crossing formulas}
\bigskip
\bigskip
\centerline{ Tian-Jun Li and Ai-Ko Liu}
\bigskip
\bigskip
\noindent{\bf  \S1. Introduction}
\medskip
The Seiberg-Witten theory was initially introduced by Seiberg and Witten
([SW], [W]) and has led to many exciting results 
 on smooth four-manifolds. When $b^+$ of a four-manifold
is larger than one, the original Seiberg-Witten
invariants constitute a map from the set of Spin$^c$ structures 
on the four-manifold to the integers. And when the first Betti number 
of the four-manifold $M$ is positive,
there is an extension of SW to a map (still denoted by SW) to
$\Lambda^*H^1(M;{\bf Z})={\bf Z}\oplus H^1\oplus\cdots\oplus\Lambda^{b_1}H^1$.
In the case when $b^+$ is equal to one, the map SW depends also
on a choice of chamber.
On a symplectic manifold, as observed in [T1], there is a canonical identification of the set of Spin$^c$
structures with $H^2(M;{\bf Z})$ and the symplectic form picks out a unique chamber.
Thus, on a symplectic manifold, SW can be viewed as a map from 
$H^2(M;{\bf Z})$ to $\Lambda^*H^1(M;{\bf Z})$. 

In a series
of remarkable works [T1-T5], Taubes developed the Seiberg-Witten theory on symplectic manifolds. 
In [T2], he defined the Gromov-Taubes invariants of
symplectic 4-manifolds counting embedded,
(but not necessarily connected) pseudo-holomorphic curves. (It was recently proved
 in [IP] that the Gromov-Taubes invariants can also be constructed
from the Ruan-Tian invariants [RT]).  
In [T3-T5] Taubes proved, on a symplectic
4-manifold with $b^+>1$, the equivalence between Seiberg-Witten invariants and 
Gromov-Taubes invariants. 
And in the case when $b^+=1$, the equivalence holds for all  classes
with pairing on all the embedded symplectic $-1$ spheres bigger than $-2$. 
  
  Taubes conjectured that the equivalence should still hold for the remaining classes
  in the case of $b^+=1$.
In this regard,   McDuff suggested a modification of the Gromov-Taubes invariants ([M]).
Adopting McDuff's modification,
we [LL2] are able to give an affirmative answer to Taubes's 
conjecture for the remaining classes using the blow up formulas for Seiberg-Witten invariants and
Gromov-Taubes invariants.

For symplectic manifolds, the equivalence between the differentiable 
invariants SW and the 
symplectic invariants GT  has led to many striking consequences to symplectic topology.
To broaden the scope of the application of the Seiberg-Witten theory, it is 
desirable to construct some kind of secondary Seiberg-Witten invariants.
It is natural to consider the parametrized Seiberg-Witten theory.
 This means, instead of a single four-manifold, we consider
 a fibre bundle with a smooth 
four-manifold as the fibre. Upon fixing a smooth family of fibrewise Spin$^c$ 
structures and a family of real self-dual two-forms,
there arises a family of Seiberg-Witten equations, and the space of
fibrewise monopoles gives rise to various kinds of invariants, the 
 family Seiberg-Witten invariants.

 This idea of  parametrized Seiberg-Witten theory was suggested by 
 Donaldson in [D2] (see also the recent papers [K] and [Rub]). 
 In this paper, we will develop
  the Seiberg-Witten theory in the family setting in some generality.
Though we restrict ourselves to the case that the base
 is a closed oriented manifold, many of the results immediately apply to
 the case where the base is an oriented manifold with boundary. 
 
 A new feature of the family Seiberg-Witten theory is
 that the chamber structure plays a more prominent role. 
For the ordinary Seiberg-Witten invariants, only if $b^+=1$ do 
they depend on chambers.
 For the family Seiberg-Witten invariants, 
they depend on chambers as long as $b^+-1$ is less than or equal
to the dimension of the base, 
 and the chamber structure is much more complicated. 
 We will introduce a finite dimensional bundle, the period bundle,
  and make use of this bundle to classify the set of chambers. 
In general the chamber structure is very complicated. However, 
  when the dimension of the base is exactly equal to $b^+-1$, the chamber structure
  is simple, it is either $\bf Z_2$ if $b^+=1$ or $\bf Z$ if $b^+>1$. 
We call this case the critical case. 

As a pleasant by-product of the study of the chamber structure, 
we can extend the scope of  a certain homomorphism $Q$ from the homology
of the space of cohomologous symplectic forms introduced by Kronheimer [K].
This simple extension is interesting; for example, on $S^2\times S^2$, 
it is nontrivial while Kronheimer's original $Q$
is not even defined. 

In the critical case we extend our techniques in [LL1] to prove a wall crossing
formula. The wall crossing formula is universal in the sense that it
only depends on the four-manifold, not on the topology of the
fibre bundle.  

But there are many interesting
families with high dimensional bases which fall into the non-critical case
$b^+-1< dim B$. We will study the Fulton-MacPherson families,
which are built up from the Fulton-MacPherson spaces and 
derive the corresponding wall crossing formulas. 
These families are intimately related to the counting of nodal  
pseudo-holomorphic curves. 
  In fact, after a certain transformation, the wall crossing formulas 
are manifestly tied to an enumeration
 problem in algebraic geometry.  This also strongly supports (as a nontrivial
 example) the general conjecture between the family Seiberg-Witten
 theory and the family Gromov theory. In fact, Ai-Ko Liu [Liu] proves a blow-up formula
 for parametrized Seiberg-Witten invariants, and when applied to 
 the Fulton-MacPherson families of K\"ahler surfaces, it is closely
 related to the wall crossing formula and has beautiful applications
 to counting curves in K\"ahler surfaces.

A  symplectic family is  
a fibre bundle with symplectic four-manifolds as fibres and a smooth family of
symplectic forms over the fibres.  
This is an extremely interesting class of families. Nice examples of 
such families are hyperk\"ahler families, or more general winding 
families (defined in section 3) of
the K3 surface and the four-torus (the bases are $S^2$).
A first observation is 
that there is a canonical chamber associated to any symplectic family,
the symplectic chamber. 
The family Seiberg-Witten invariants in the symplectic chamber will be called
the family Taubes-Seiberg-Witten invariants. 
 Though the Seiberg-Witten
invariants for the K3 surface and the four-torus are trivial, the invariants
of the winding family are very rich. 
In fact we realized that the counting of curves on a projective
K3 surface is closely related to the invariants (Seiberg-Witten and Gromov)
for the hyperk\"ahler family, and Yau and Zaslow's beautiful conjectural
generating function of curves on projective K3 [YZ]
was one of the motivations to develop the parametrized
Seiberg-Witten theory. Recently, Yau and Zaslow's formula has been indeed
confirmed via calculation of the Gromov invariants
of a hyperk\"ahler family of K3 surfaces by Bryan and Leung [BL1].

Because of Taubes's results [T1-T5], one naturally expects that
the family Seiberg-Witten theory is particularly interesting for symplectic 
families and the nonvanishing of the family Taubes-Seiberg-Witten invariants will imply
the existence of fibrewise pseudo-holomorphic curves. 
This is indeed the case.
However,  in setting  up the family Gromov-Taubes
invariants counting fibrewise pseudo-holomorphic curves and pushing
Taubes's equivalence between SW and GT to the family setting, 
we have had only  limited success.
For a restricted class of families, including the winding families of
 the four-torus and primary Kodaira surfaces, 
we can give a good definition and prove TSW=GT.
These results will appear in [LL3] and be applied to study symplectic manifolds
with torsion symplectic canonical classes.

The organization of the paper is as follows. In 
\S2, we introduce the family Seiberg-Witten invariants.
In \S3, we introduce the period bundle and study the chamber structure.
    In \S4, we
 prove the wall crossing formula in the critical case and study
some interesting examples.  
We also derive the wall crossing formula for the Fulton-MacPherson family.

The authors wish to thank  S. T. Yau who suggested to
us to study parametrized Seiberg-Witten theory
soon after the birth of the Seiberg-Witten theory
and has given us many valuable ideas and shown great interest.
The authors wish to thank C. Taubes and  G. Tian 
for many suggestions and sharing their ideas. 
 The authors wish to thank  E. Brown, R. Fintushel,  G. Moore,  R.
Lee,  G. Zuckerman for their interest in this work. Both authors are
partially supported by NSF and wish to acknowledge the support and hospitality of 
IAS. Finally the authors wish to thank  A. Greenspoon and the referee for 
careful reading and excellent suggestions which greatly improve the presentation of
the paper.

\bigskip
\bigskip
\noindent{\bf \S2. Family Seiberg-Witten invariants}
\medskip
In this section, we set up the family Seiberg-Witten theory. 
In the first two subsections, the family Seiberg-Witten equations and
the family Seiberg-Witten invariants are introduced.  The last 
subsection is devoted to symplectic families.
\medskip
\noindent{\bf \S2.1. Family Seiberg-Witten equations}
\medskip
Let $M$ be an oriented closed 4-manifold. 
Let $b_i$ denote the $i-$th Betti number of $M$ and
let $b^+$ denote the dimension of a maximal subspace, $H^+(M;{\bf R})\subset
H^2(M;{\bf R})$ on which  the cup product form is positive. 

Let $B$ be an oriented 
closed manifold and ${\cal X}$ be the total space of a fibre bundle
with fibre $M$ and base $B$.  

Denote the tangent bundle along the fibres
by $T({\cal X}/B)$ and the bundles of $i-$forms along the fibres
by ${\Lambda^i}$. Choose a metric $G$ on $T({\cal X}/B)$;
then it defines a principal $SO(4)$ bundle of frames, ${\cal F}r\rightarrow 
{\cal X}$ 
and two associated bundles,
 $\Lambda^{+}$ of self-dual 2-forms along the
fibres  and ${\Lambda}^{-}$ of anti self-dual 2-forms along the fibres.

For any $b\in B$, denote by ${\cal X}_b$ the fibre at $b$.
It will be understood that for any object defined
on the total space ${\cal X}$, the object with subscript $b$ will denote
the restriction to the fibre ${\cal X}_b$.

A Spin$^c$ structure ${\cal L}$ on ${\cal X}/B$ is an equivalence
class of lifts of ${\cal F}r$ to a principal
Spin$^c(4)$ bundle ${\cal F}$. Recall that the group Spin$^c(4)$ is the group 
$(SU(2)\times SU(2)\times U(1))/\{\pm1\}$ and the group $SO(4)$ is also the same as
the group $(SU(2)\times SU(2))/\{\pm 1\}$, and the homomorphism 
$Spin^c(4)\longrightarrow SO(4)$ is simply the map forgetting the factor
$U(1)$.

Associated with ${\cal F}$, there are two $U(2)$ bundles
of spinors ${\cal S}^+$ and ${\cal S}^-$ coming from 
the two natural homomorphisms from Spin$^c(4)$ to $U(2)$ which is
the same as $(SU(2)\times
U(1))/\{\pm1\}$. ${\cal S}^+$ is distinguished by the
identification of its projective bundle with the unit two-sphere bundle of 
${\Lambda}^{+}$. The natural homomorphism $ T^*({\cal X}/B)\otimes {\cal S}^+
\longrightarrow {\cal S}^-$ defines the Clifford multiplication.
And the adjoint of the Clifford multiplication endmorphism defines
  a canonical homomorphism
$$\tau:\hbox{End}({\cal S}^+)\longrightarrow \Lambda^+\otimes {\bf C}.$$

Let $L=det({\cal S}^+)$. $L$ is viewed as a family of $U(1)$ bundles
over $B$. Given a smooth family of
$U(1)$ connections $A=\{A_b\}$ on $L$,
combined with the family of Levi-Civita connections on ${\cal F}r$, 
it defines a smooth family of covariant derivatives on ${\cal S}^{+}$, still denoted by
$A=\{A_b\}$. Let $\psi$ be a section of ${\cal S}^+$; then it is naturally viewed
as a smooth family of sections of ${\cal S}^+_b$.
Denote by ${\cal C}_b$ the space of pairs $(A_b,\psi_b)$
where $A_b$ is a $U(1)-$connection on $L_b$ and $\psi_b$ is a section 
of ${\cal S}^+_b$. 
${\cal C}_b$ is an affine Fr\'echet manifold modeled on $i\Lambda_b^1\times
C^{\infty}({\cal S}_b^+)$ ($\Lambda_b^1$ is the space of smooth 1-forms
on ${\cal X}_b$).
Then the configuration space ${\cal C}$ is the space of pairs $(A,\psi)$.
 It is just the space of sections of the
infinite dimensional affine bundle over $B$ with fibre ${\cal C}_b$.
Call an element $(A,\psi)\in {\cal C}$ irreducible if $\psi_b$ is not identically
zero for any $b$. Denote by ${\cal C}^*$ the irreducible part of ${\cal C}$.

The unperturbed Seiberg-Witten equations on ${\cal X}/B$ are equations for
$(A,\psi)\in {\cal C}$:
$$\eqalign{D_{A_b}\psi_b&=0\cr
P^+_bF_{A_b}&={1\over 4}\tau(\psi_b\otimes \psi_b^*).\cr}\eqno (2.1)$$
In the first equation, $D_{A_b}:\Gamma({\cal S}^+_b)\longrightarrow \Gamma({\cal S}^-_b)$
is the Dirac operator, a first order differential operator
defined using the Clifford multiplication and the covariant 
derivatives $A_b$ on ${\cal S}^{+}_b$. 
In the second equation, $P^+_b:\Lambda^2_b\longrightarrow \Lambda^+_b$
is the orthogonal projection and $F_{A_b}$ is the curvature two-form of $A_b$.

Let $\mu\in \Lambda^+$, then $\mu_b$ defines a real valued
self-dual two-form on ${\cal X}_b$. It is quite useful to consider 
perturbations which are of the form 
$$\eqalign {D_{A_b}\psi_b&=0\cr
P^+_bF_{A_b}&={1\over 4}\tau(\psi_b\otimes \psi_b^*)+i\mu_b.\cr}\eqno (2.2)$$

The parameter space for the family Seiberg-Witten equations is thus  the following
subspace of the product
${\cal I}\times i\Lambda^2$:
$$\Gamma=\{(G, \mu)|\hbox{ for each $b$, } 
\mu_b \hbox{ is self-dual with respect to $G_b$}\}$$
where ${\cal I}$ is the space of metrics on $T({\cal X}/B)$.
${\cal I}$ is an infinite dimensional fibre bundle over $B$
whose fibre at $b$ is the space of Riemannian metrics on ${\cal X}_b$
(a Fr\'echet manifold).

$\Gamma$ has an alternative description which is often more convenient. 
Denote by $\tilde{\cal P}$ the infinite dimensional bundle
over $B$ whose fibre over $b$ is the space of pairs $(G_b, \mu_b)$,
where $G_b$ is a metric on ${\cal X}_b$ and $\mu_b$ is a real valued two-form and self-dual with respect to 
$G_b$. Then $\Gamma$ is just the space of
sections of $\tilde{\cal P}$.
To study $\Gamma$, we introduce $\bar{\cal P}$, a finite dimensional 
bundle over $B$ whose fibre 
over $b$ is the product of the real Grassmannian
$Gr(H_b, b^+)$  and the vector space $H_b$, where $H_b$ is 
$H^2({\cal X}_b;{\bf\ R})$. There is a natural map $pr_{\cal L}$ from 
$\tilde{\cal P}$ to $\bar{\cal P}$
sending each pair $(G_b, \mu_b)$ to the pair 
$([{\cal H}^+], {\cal H}(2\pi c_1(L)-\mu_b))$
where $[{\cal H}^+]$ is the point in the real Grassmannian given by the
$b^+-$dimensional subspace of self-dual harmonic two-forms and
${\cal H}(2\pi c_1(L)-\mu_b)$ is the harmonic projection of 
$2\pi c_1(L)-\mu_b$, both with respect to
the metric $G_b$.
The next section will be devoted to the homotopy classes of the sections
of the period bundle ${\cal P}$.

 Given any pair $(G, \mu) \in \Gamma$, denote the space of solutions 
to (2.2) by 
${\tilde{\cal M}}$ and its irreducible
part by $\tilde{\cal M}^*$. 

The gauge group is replaced by the bundle of
groups ${\cal G}$ over $B$, whose fibre over $b$ is 
the group ${\cal G}_b=C^{\infty}({\cal X}_b, U(1))$.
${\cal G}_b$ acts smoothly on ${\cal C}_b$ by sending a 
map $g_b$ and a pair $(A_b, \psi_b)$ to $(A_b+2g_bdg_b^{-1}, g_b\psi_b)$
 and the action is free  on ${\cal C}_b^*$.
This action naturally extends to a smooth 
map from ${\cal G}\times_B {\cal C}$ to 
${\cal C}$ which defines a ${\cal G}$ action on ${\cal C}$ (we abuse
language here and elsewhere in this paper; though
 ${\cal G}$ is not really a group, we will say that it acts on ${\cal C}$). 
And the space of the orbits of this action, denoted by ${\cal C}/{\cal G}$,
is given the quotient topology.
Just as in the case of the ordinary Seiberg-Witten
equations, ${\cal G}$ preserves $\tilde{\cal M}$ and 
$\tilde{\cal M}^*$, and the quotients are called ${\cal M}$ and ${\cal M}^*$
respectively with the subspace topology inherited from
${\cal C}/{\cal G}$.

For any $(G,\mu)\in \Gamma$, ${\cal M}$ is always compact.
Since $B$ is compact, the Weitzenb\"ock formula still gives a
uniform bound on the $C^0$ norm of $\psi_b$ of any solution
$(A_b,\psi_b)$ to (2.2). From it, we can derive the uniform bound on
the $C^{\infty}$ norm of solutions to (2.2) up to gauge, and the compactness
follows.

Let $(A_b, \psi_b)$ be a solution to (2.2). If we linearize the SW equations
at $(A_b, \psi_b)$ and take into account the gauge group action, we get a linear Fredholm operator,
$${\cal D}: T_bB\oplus i\Lambda^1_b\oplus \Gamma({\cal S}^+_b)
\longrightarrow i\Lambda^+_b\oplus \Gamma({\cal S}^-_b)\oplus i\Lambda_b^0.$$
${\cal D}$ is Fredholm of index
$$d=d({\cal L})=dim B+{1\over 4}[c_1(L)^2-(2e(M)+3\sigma(M))]\eqno (2.3)$$
where $e(M)$ is the Euler characteristic of $M$ and $\sigma(M)$ is the
signature of $M$.

By the Smale-Sard theorem, 
for generic $(G,\mu)\in \Gamma$ (here generic means a Baire set),
${\cal D}$ is surjective at all the irreducible solutions,
so  ${{\cal M}^*}$ is a smooth manifold of dimension $d$ as in (2.3).

Given a Spin$^c$ structure ${\cal L}$, we define a subset ${\cal WA}({\cal L})$ in ${\Gamma}$ which we will 
call the ``${\cal L}-$wall'' in 
 ${\Gamma}$. A pair $(G, \mu)$ is in ${\cal WA}({\cal L})$ if for some $b$, 
$2\pi c_1(L_b)-\mu_b$ is anti-self-dual with respect to $G_b$.  
For any pair 
 $(G,\mu)\in \Gamma-{\cal WA}({\cal L})$, $\tilde{\cal M}$
contains no reducible solutions, so $\tilde{\cal M}=\tilde{{\cal M}^*}$.

From the second description of $\Gamma$, it is not hard to see that
${\cal WA}({\cal L})$ is a codimension $max(b^+-dimB, 0)$ submanifold in $\Gamma$.
In particular, if $b^+-1\geq
dim B$, for a generic pair $(G, \mu)$ (here generic means an
open set off a submanifold of positive codimension), there are no reducible solutions.
Even when $b^+\leq dim B$,  it is not hard to 
see that if there is a pair $(G, \mu)$ such that $\mu_b$ is not exact
for each $b$, then  
any 
pair $(G, t\mu)$ with sufficient large $t$ does not lie in 
${\cal WA}({\cal L})$, thus $\Gamma-{\cal WA}({\cal L})$ still contains a
nonempty open set.
Call each connected component of $\Gamma-{\cal WA}({\cal L})$ an
${\cal L}$-chamber.

A choice of  an orientation of
the real line $det^+=det(H^0(M;{\bf R})\otimes H^1(M;{\bf R})\otimes H^+(M;{\bf R}))
\otimes H^{dimB}(B;{\bf R})$ serves to orient ${\cal M}$.

Hence there are always pairs $(G,\mu)\in \Gamma$ such that 
${\cal M}={\cal M}^*$ and ${\cal M}^*$ is a closed 
manifold of dimension $d$ as in (2.3) in the space ${\cal C}^*/{\cal G}$,
 and after fixing an orientation of the
determinant line $det^+$, it represents
 a homology class 
$[{\cal M}^*]$ in
$H_d({\cal C}^*/{\cal G};{\bf Z})$.
If $b^+-1\geq dim B$, such pairs $(G,\mu)\in \Gamma$ are generic.

\bigskip
\noindent{\bf \S2.2 Family Seiberg-Witten invariants}
\medskip

Denote the set of Spin$^c$ structures on ${\cal X}/B$ 
by ${\cal SP}$. Although the definition of ${\cal SP}$ requires a
choice of metric on $T({\cal X}/B)$, there is a natural identification
between such sets defined by any two metrics. Thus, the set ${\cal SP}$ only
depends on the fibre bundle ${\cal X}/B$ and is an affine space 
modeled on $H^2({\cal X};{\bf Z})$.

Given an element ${\cal L}\in {\cal SP}$, the definition of the family SW
invariants of ${\cal L}$ depends on a cohomology class $\Theta$ on the space
${\cal C}^*/{\cal G}$
and a choice of parameter $(G,\mu)\in \Gamma$.

Fix a pair $(G, \mu)$ such that ${\cal M}$ is a closed manifold of 
dimension $d=d({\cal L})$ as in (2.3). Fix an orientation of the determinant line.
Choose $\Theta\in H^{l}({\cal C}^*/{\cal G};{\bf Z})$; the Seiberg-Witten invariant associated to $\Theta$ is 
defined to be
$$SW({\cal X}/B, {\cal L}, \Theta)=<\Theta, [{\cal M}]>\eqno (2.4)$$
if $l=d$, and zero otherwise, where $<-,->$ is the pairing between
cohomology and homology on the space ${\cal C}^*/{\cal G}$.
Sometimes we will drop ${\cal X}/B$ or ${\cal L}$ in the notation when there is 
no confusion.

If $b^+-1>dimB$, given a Spin$^c$ structure ${\cal L}$ on ${\cal X}/B$
then the SW invariants of ${\cal L}$ are defined for generic 
pair $(G,\mu) \in \Gamma -{\cal WA}({\cal L})$. And since
$\Gamma-{\cal WA}({\cal L})$ is path connected, from the
 standard  cobordism arguments, 
the value of the Seiberg-Witten invariant is independent of the choice of
the parameter. This is analogous to the invariance of the
ordinary Seiberg-Witten invariant in the case $b^+$ is greater than one.
Hence we have the following 
theorem:

\proclaim Theorem 2.1. Let $M$ be a closed oriented 4-manifold and $B$ be
a closed oriented manifold. Let ${\cal X}$ be a fibre bundle with
fibre $M$ and base $B$, ${\cal L}$ be a Spin$^c$ structure on
${\cal X}/B$ and $\Theta$ be a cohomology class in
$H^*({\cal C}^*/{\cal G};{\bf Z})$. Fix an orientation of $det^+$.
If $b^+-1>dim B$, then $SW ({\cal X}/B, {\cal L}, \Theta)$
is independent of the choice of generic $(G,\mu)$, and hence is a 
differentiable invariant of the fibre bundle. Furthermore, if $f$ is a
self-diffeomorphism of ${\cal X}$
preserving the fibres, then 
$$SW({\cal X}/B, {\cal L}, \Theta)=\pm SW({\cal X}/B, f^*{\cal L},
f^*\Theta).$$\par

In the case $b^+-1\leq dim B$, the preceding theorem does not hold.
However, by the standard arguments, the  following conclusion still holds:

\proclaim  Theorem 2.2. 
 Let $M$ be a closed oriented 4-manifold and $B$
a closed oriented manifold. Let ${\cal X}$ be a fibre bundle with
fibre $M$ and base $B$, ${\cal L}$ be a Spin$^c$ structure on
${\cal X}/B$ and $\Theta$ be a cohomology class in $H^*({\cal C}/{\cal G}; {\bf Z})$.
 Fix an orientation of $det^+$.
Suppose $b^+-1\leq dim B$ and let $\Gamma_c$ be an ${\cal L}-$chamber. 
 Then  if $0< b^+\leq dimB +1$,  
  $SW ({\cal X}/B, {\cal L}, \Theta, c)$ is 
defined for generic pairs $(G, \mu)\in \Gamma_c$ and is independent of
the choice of $(G, \mu)$.\par

With Theorem 2.2 understood, it is important to pin down the
 dependence on chambers. This will be the content of the next two sections.

In the remaining part of this subsection, we will describe how
to construct some interesting cohomology classes on
the infinite dimensional object ${\cal C}^*/{\cal G}$.

The first observation is that it is 
homotopic to a $CP^{\infty}\times T^{b_1}$  bundle
over $B$. 
Any cohomology class of $B$ pulls back to a cohomology class
of ${\cal C}^*/{\cal G}$.
If $\Theta$ can be factored by such a class with degree $l$, then we say
$\Theta$ has weight $l$.
There is a distinguished class $[B]^{\wedge}$ in 
$H^{\dim B}({\cal C}^*/{\cal G};{\bf Z} )$
 which is the pull back of the fundamental cohomology class of $B$.
If $\Theta$ is factored by $[B]^{\wedge}$,  then the family Seiberg-Witten invariant
becomes an ordinary Seiberg-Witten invariant. Similarly, 
if $\Theta$ has weight $l$, then the invariant reduce to an invariant
on a codimension $l$ family.

Another observation is that
${\cal G}$ acts naturally  on the complex line bundle $L$ over ${\cal X}$.
The diagonal action on the product $L\times {\cal C}^*$ defines a
complex line bundle $L\times_{\cal G}{\cal C}^*$ over
 ${\cal X}\times ({\cal C}^*/{\cal G})$.
Denote the first Chern class of this line bundle by $u$ and the restriction of $u$ to
 ${\cal C}^*/{\cal G}$ by $H$. $H$ will be called
the hyperplane class. 

Via slant product with $u$, we get a map 
$$\mu: H_i({\cal X}; {\bf Z})\
\longrightarrow H^{2-i}({\cal C}^*/{\cal G};{\bf Z}).$$ 
Pick homology classes $\zeta_1, \cdots, \zeta_l$ in $H_{*}({\cal X};{\bf Z})$;
we can define 
$$ SW({\cal X}/B, {\cal L}, \mu(\zeta_1)\cdots
\mu(\zeta_l))=SW({\cal X}/B, {\cal L}, \mu(\zeta_1)\wedge\cdots\wedge
\mu(\zeta_l)).$$

A particular important invariant, the pure invariant, is given in the 
following definition.

\proclaim Definition 2.3 (pure invariant). 
When the dimension of ${\cal M}$ is even, we
define the pure invariant to be 
$$SW(H)=\int_{\cal M} H^{dim {\cal M}/2}.$$

 Set $U=H\cup  [B]^{\wedge}$. 
Via slant product with $U$, we get a map 
$$\bar \mu: H_i({\cal X}; {\bf Z})\
\longrightarrow H^{dimB+2-i}({\cal C}^*/{\cal G};{\bf Z}).$$ 
Pick homology classes $\xi_1, \cdots, \xi_l$ in $H_{*}({\cal X};{\bf Z})$;
we can define 
$$SW({\cal X}/B, {\cal L} , \bar \mu(\xi_1),\cdots, \bar \mu(\xi_l)).$$
It counts the number of fibrewise monopoles whose zero set intersect
each of the submanifolds $V_i$, where $V_i$ represents $\xi_i$.

When ${\cal X}$ is a product bundle, $\mu$ and $\bar \mu$ are simply related by
$$\mu(\zeta)=\bar \mu (\zeta\times [B]).$$

\noindent{\bf Example 2.4} (section invariants and circle section invariants). Suppose ${\cal X}$
has sections. Given any section $s:B\longrightarrow {\cal X}$, we can define
a based gauge group ${\cal G}(s)\subset {\cal G}$ consisting of
maps sending $s(B)$ to the identity. Define
$${\cal C}(s)={{\cal C}}/{\cal G}(s)\hbox{\quad
 and\quad}{\cal C}^*(s)={{\cal C}^*}/{\cal G}(s).
$$
Then ${\cal C}^*(s)$ is a principal $S^1$ bundle over ${\cal C}^*/{\cal G}$
and its Euler class $e(s)$ is in $H^2({\cal C}^*/{\cal G};{\bf Z})$.
Notice that a section $s$ determines a homology class in $H_*({\cal X};{\bf Z})$,
and $e(s)$ is simply the slant product of this homology class with $U$.
If $d$ is an even integer, then after picking $d/2$ 
sections, $s_1, \cdots, s_{d/2}$, the cohomology class 
$\prod_{i=1}^{d/2}e(s_i)$ is of degree $d$.
We can define the corresponding Seiberg-Witten invariant $SW({\cal X}/B, {\cal L}, s_1,s_2,\cdots, s_{d/2}$).
 This invariant only
depends on $[s_i]$, $i=1, \cdots, d/2$, the homotopy classes of $s_i$,
as the corresponding cohomology classes depend only on the homotopy
classes. 
A geometric interpretation of this invariant is that it
counts the number of fibrewise monopoles vanishing along these sections.

Similarly, if ${\cal X}$ allows sections of circles, 
 each circle section defines a homology class of ${\cal X}$,
and induces a cohomology class of ${\cal C}^*/{\cal G}$
via the slant product. Thus, given $d$ circle sections
$\gamma_1,\cdots, \gamma_d$, we can define the SW invariant
$SW({\cal X}/B, {\cal L}, \gamma_1,\cdots, \gamma_d)$
which  only depends on the 
isotopy classes of the circle sections. 
Geometrically, this invariant counts the number of fibrewise monopoles
which are not nowhere vanishing along each of the circle sections.
More generally, given an integer $p$
such that $d-p$ is a nonnegative even integer, and $p$ circle sections 
$\gamma_1,\cdots,\gamma_p$ and $(d-p)/2$ sections
$s_1,\cdots, s_{(d-p)/2}$, we can define the Seiberg-Witten invariant
$SW({\cal X}/B, {\cal L}, \gamma_1,\cdots, \gamma_p,s_1,\cdots, s_{(d-p)/2})$.

\noindent{\bf Remark 2.5}. Donaldson, in his beautiful survey article [D2],
suggested the possibility of constructing the family Seiberg-Witten 
invariants. He further suggested that these invariants should be
viewed as cohomology classes of $Bdiff(M)$, the classifying space
of the group of diffeomorphisms of $M$.

\bigskip
\noindent{\bf \S2.3. Seiberg-Witten invariants on symplectic families}
\medskip
In this subsection, we describe the special features of 
the family Seiberg-Witten invariants on symplectic families.
The readers should consult [T1].

Let ${\cal X}$ be a fibre bundle of a  four-manifold
over $B$ and $\omega$ be a two-form on $\cal X$ whose restriction
to each fibre is a symplectic form. Each fibre ${\cal X}_b$ is oriented by
$\omega_b\wedge \omega_b$.

Fix a metric $G$ such that $\omega_b$ is self-dual with
respect to $G_b$. $G$ can be further normalized such that
$\omega_b$ has length $\sqrt 2$. Such a family metric
is called an $\omega$ compatible metric. $\omega$ and 
an $\omega$ compatible metric determine
a smooth family of compatible almost complex structures $J$ on
the fibres, hence
a complex structure (still denoted by $J$)
of the bundle $T({\cal X}/B)$. Taking the complex determinant,
one obtains a
 complex line bundle $K$ over ${\cal X}$,
whose restriction to each
fibre ${\cal X}_b$ is given by det($T^{0,1}{\cal X}_b)$.

There is a canonical Spin$^c$ structure with the associated
bundles ${\cal S}^+$ and ${\cal S}^-$ naturally isomorphic to
$\bracevert\oplus K^{-1}$ and
$\Lambda^{0,1}=\Lambda^{0,1}(T({\cal X}/B))$, where $\bracevert$ is the trivial
complex bundle over ${\cal X}$. The splitting of ${\cal S}^+$ is induced by the
Clifford action of $\omega$, which has eigenvalue $-2i$ on the trivial
summand $\bracevert$ and eigenvalue $2i$ on the $K^{-1}$
summand.

This canonical Spin$^c$ structure induces a natural identification
between the set ${\cal SP}$ and $H^2({\cal X};{\bf Z})$.
Under this identification, a class $e\in H^2({\cal X};{\bf Z})$
is sent to the Spin$^c$ structure whose
${\cal S}^{\pm}$ bundles are given by
$${\cal S}^+=E\oplus (K^{-1}\otimes E)\hbox{\quad  and \quad }
{\cal S}^-=\Lambda^{0,1}\otimes E.\eqno (2.5)$$
where $E$ is a complex line bundle over ${\cal X}$ whose first Chern class is
$e$.

There is a natural orientation for the line $det^+$ provided
an orientation of $B$ is fixed. This is because the half-self-dual
complex on $M$ is within a relatively compact perturbation of a
complex linear complex arising from a compatible almost complex structure.

To introduce the family version of 
Taubes's perturbation of the SW equations it
requires the introduction of a smooth family of canonical
connections $A^0$ (up to the gauge action) on $K^{-1}_b$.
A family of connections on $K^{-1}$, coupled with the family of
Levi-Civita connections, gives rise to a family of covariant
derivatives, $\nabla_A$, on ${\cal S}^+$. Through restriction and projection,
$\nabla_A$ defines a family of covariant derivatives, 
$\nabla_A$, on the trivial
summand $\bracevert$.
The family $A^0$ is characterized by the property that the corresponding
family of covariant derivatives, $\nabla_{A^0}$, admits a family of
non-trivial covariantly constant sections, $u^0$. In the subsequent
discussions,  $u^0$ will be normalized such
that $u^0_b$ has norm one.

Let ${\cal L}$ be a Spin$^c$ structure on
${\cal X}/B$ specified by a complex line bundle
$E$ over ${\cal X}$ under the aforementioned identification.
Recall, in the family setting, that the SW equations
are equations for pairs $(A_b, \psi_b)$,
where $A_b$ is a connection on the complex line bundle $L_b=det {\cal S}^+_b$
and $\psi_b$ is a section of ${\cal S}^+_b$. Note that each
line bundle $L_b$
is naturally the restriction of the complex line bundle $L=det{\cal S}^+$
to the fibre ${\cal X}_b$. From (2.5), $L_b$ is $K^{-1}_b\otimes E_b$.
Thus, with $A^0$ fixed, a connection $A_b$ on $L_b$ is
written as $A_b=A^0_b+2a_b$, 
where $a_b$ is a connection on the complex line bundle $E_b$.

Following Taubes, we choose the one-parameter family
(parametrized by $r$) of the perturbation two-forms $\mu$ by
$$\mu_b={-r\over 4}\omega-i P_+F_{A^0_b}.\eqno (2.6)$$

It is useful to rewrite the spinor corresponding to the splitting 
and the parameter $r$ as
$$\psi_b=r^{1/2}\cdot (\alpha_b u^0_b+\beta_b).$$

Now the equations read
$$\sigma(u^0\otimes \nabla_{a_b}\alpha_b)+D_{A_b}\beta_b=0,$$
$$P_+F_a=-{i\over r}\cdot (1-|\alpha_b|^2+|\beta_b|^2)\cdot \omega_b
+{ir\over 4}(\alpha_b \beta_b^* +\alpha_b^*\beta_b).\eqno (2.7)$$
Here, $\alpha_b\beta_b^*$ and $\alpha_b^*\beta_b$, being respectively sections of
$K_b$ and $K_b^{-1}$, are identified as 
sections of $\Lambda_b^+\otimes {\bf C}$.

Notice that there is a canonical chamber, called the Taubes chamber,
for any Spin$^c$ structure ${\cal L}$ on a symplectic family.
On each fibre ${\cal X}_b$, when $r$ is sufficiently large,
$2\pi c_1(L)- \mu_b$ has positive square due to the term
$r^2\omega\wedge\omega/16$. So no reducible solutions
can possibly occur for all large $r$. Since $B$ is compact,
there is a uniform constant such that $2\pi c_1(L)-\mu_b$
has positive square for all $b$ once $r$ is greater than
that constant.   This means that for all large $r$,
the pairs $(G, \mu)$ with $\mu$ given by (2.6)  lie in the
same ${\cal L}$ chamber, and this chamber will be
called the Taubes chamber.
Notice that this chamber only depends on the family of symplectic forms
through its deformation class: If $\omega$ and $\omega'$
are two families of symplectic forms on a fibre bundle, and
there is a path of families of symplectic two-forms
connecting $\omega$ and $\omega'$, then the Taubes chambers for
these two families are the same. Finally, let us give a definition.

\proclaim Definition 2.6. Given a symplectic family, the Seiberg-Witten
invariants in the Taubes chamber are called the Taubes-Seiberg-Witten
invariants.

\bigskip
\noindent {\bf \S3. The period bundle and the chamber structure}
\medskip

When $0\leq b^+-1\leq dimB$, $\Gamma-{\cal WA}({\cal L})$ may have more than one
component, i.e. the number of ${\cal L}$ chambers may be greater than one.
 When $B$ is a point, we
are back to  the ordinary Seiberg-Witten theory
in the case $b^+=1$, and the chamber structure is
well understood; there are two ${\cal L}$ chambers.
For general $B$, the chamber structure is much more complicated.
 In this section we will discuss the chamber structure 
when $0\leq b^+-1\leq dim B$ and 
give a classification in cohomological terms.
We will derive the
wall crossing formula in the next section. 

In the first subsection, we will reduce this 
problem to a finite dimensional bundle 
${\cal P}$, the period bundle, and study the topology of
this bundle.
In the next subsection, we will give the classification. 
\bigskip
\bigskip
\noindent{\bf \S3.1. Reduction to the period bundle ${\cal P}$}
\medskip

The parameter space $\Gamma$ is the space of sections of the
infinite dimensional bundle $\tilde{\cal P}$. 
Fix a Spin$^c$ structure ${\cal L}$ with the associated complex 
 line bundle $L$; there is a natural map
$pr_{\cal L}$ from $\tilde {\cal P}$ to $\bar{\cal P}$  sending the pair
$(G_b,\mu_b)$ to the pair 
$([{\cal H}_b^+], {\cal H}(2\pi c_1(L_b)-\mu_b))$,
where $[{\cal H}_b^+]$ is the point in the real Grassmannian given by the
$b^+-$dimensional subspace of self-dual 
harmonic two-forms and
${\cal H}_b(2\pi c_1(L_b)-\mu_b)$ is the harmonic projection of 
$2\pi c_1(L_b)-\mu_b$, both with respect to
the metric $G_b$.

The map $pr_{\cal L}$ is not surjective; it surjects onto an open subbundle.
In that regard, we introduce the following definition. 

\proclaim Definition 3.1. (self-dual Grassmannian). Let $V$ be a vector space with
 a nondegenerate quadratic form defined over $\bf R$ with signature $(p,q)$.
 All the positive definite $p-$planes in $V$
 form an open subset of the $p-$plane Grassmannian in $V$.
It will be called the self-dual Grassmannian, and denoted by
$G_p^+(V)$.

Let ${\cal P}^+$ denote the subbundle of $\bar{\cal P}$ whose fibre 
over $b$ is 
$G_{b^+}^+(H^2({\cal X}_b;{\bf R}))\times H^2({\cal X}_b;{\bf R})$.
We first give a simple observation.

\proclaim Lemma 3.2.
The map $pr_{\cal L}$ surjects onto ${\cal P}^+$ and is regular everywhere.
Moreover, the preimage of every point is contractible.

These properties are not hard to prove. We just mention two useful facts here.
 The first fact is that the space of metrics fibres over the space of conformal structures with
contractible fibres, the space of positive functions. 
 The second fact is that  a conformal structure is specified by a real three-dimensional 
positive subbundle, and relative to a fixed Riemannian metric,
the space of conformal structures is identified with 
bundle maps 
$$\nu:\Lambda^+\longrightarrow \Lambda^-, |\nu(\eta)|<|\eta|$$

The point is now that the study of the parameter space $\Gamma$,
which is the space of  sections of $\tilde{\cal P}$, is reduced to the
study of the space of sections of the finite dimensional 
bundle ${\cal P}^+$.


Given a Spin$^c$ structure ${\cal L}$, 
what we really want to 
 understand is the space $\Gamma-{\cal WA}({\cal L})$.
For this purpose, we need to introduce
the period bundle, which is a subbundle of
${\cal P}^+$.

Let us first assume that $B$ is a point. 
 Then
 a pair $(G,\mu)$ lies in ${\cal WA}({\cal L})$ if
 $P_+F_{A}-i\mu=0$ has a solution. 
But this is equivalent to the following: inside $C^{\infty}
({\Lambda}^+)$, the affine subspace 
$\{A\longrightarrow {1\over i}F_A-\mu\}$ passes through the 
origin,  and thus coincides with the linear subspace 
$dC^{\infty}(\Lambda^1)$.
This condition is also equivalent to 
${\cal H}(2\pi c_1(L)-\mu)$ being perpendicular to ${\cal H}^+$. 
This simple observation enables us to describe explicitly the
image of $\Gamma-{\cal WA}({\cal L})$ under the map $pr_{\cal L}$.

As $B$ is a point, $\tilde {\cal P}$ is 
itself the space $\Gamma$ and ${\cal P}^+$ is simply 
the  product $G_p^+(V)\times V$ where $V=H^2(M;{\bf R})$. 
 The question we are interested in is what kind of positive $p-$plan will be
perpendicular to a chosen element $l$ in the vector space $V$.  
In that regard, we introduce a subspace of ${\cal P}^+$,
 ${\cal P}$, consisting of the pairs $(W,l)$ in ${\cal P}^+$
such that $l$ is not perpendicular to $W$. ${\cal P}$
 is called the
period space of $V$.
Under the map $pr_{\cal L}$, $\Gamma-{\cal WA}({\cal L})$
and  ${\cal WA}({\cal L})$ are sent onto ${\cal P}$ and ${\cal P}^+
-{\cal P}$ respectively. As remarked before, the
map $pr_{\cal L}$ is regular everywhere with 
contractible preimages, therefore $\Gamma-{\cal WA}({\cal L})$
is homotopic to ${\cal P}$. 

To study the topology of ${\cal P}$, it is necessary to first
understand $G_p^+(V)$.
If $p$ is nonzero, $G_p^+(V)$ is actually contractible; this follows from the 
the fact that
$G_p^+(V)$ can be identified with $SO(p,q)/SO(p)\times SO(q)$ which is a
homogeneous space formed by a maximal compact subgroup quotient. Another
way to check it is by observing that the projection from one positive $p-$plane to another one  induces an 
isomorphism between these two vector spaces.
Fixing a base positive $p-$plane and projecting every
 other positive $p-$plane to the fixed one induces a smooth
function on $G_p^+(V)$ which has a unique maximum at the base point. This can
be done by comparing the 
volume form (absolute value) of a positive $p-$plane
 with that of the base point. 
The upward gradient flow associated to this function
produces a homotopy which
 shrinks $G_p^+(V)$ to a point.

The period space, viewed as an open submanifold in $G_p^+(V)\times
 V$, can be characterized as the complement of the vanishing locus 
of certain explicit map from $G_p^+(V)\times V$ to ${\bf R}^p$. 
The map is defined by first 
choosing a base point in $G_p^+(V)$, then projecting the specified 
$p-$plane 
  to the other positive $p-$plane. As we have explained, the projection
 map is always nondegenerate; thus a frame in the base $p-$plane will induce 
a smooth family of frames.  Now we are ready to define the map
from $G_p^+(V)\times V$ to ${\bf R}^p$. Picking a point in $G_p^+(V)$,
which is just a positive $p-$plane,  and an element in $V$,
  the pairing between the attached frame in the $p-$plane
 and the element in $V$
 gives us $p$
 elements in ${\bf R}$; thus it determines a map $c$ from $G_p^+(V)\times V$
 to ${\bf R}^p$. The pairs $(W, l)$ in the preimage of $0\in {\bf R}^p$ are
 characterized by the property 
that the $p-$plane $W$ is perpendicular to $l$ in $V$. Deleting
 these points from $G_p^+(V)\times V$ results in the period space.
This description of the period space readily leads to
the following proposition. 

\proclaim Proposition 3.3. Let $V$ be a vector space with a quadratic form
 of signature $(p,q)$ defined over ${\bf R}$. The period space associated to the 
vector space $V$ is homotopic equivalent to $S^{p-1}$.

In fact it is easy to show that the map 
$c$ is regular everywhere. Thus the preimages of all point
 in ${\bf R}^p$ are diffeomorphic to each other. Therefore, the period space 
has a structure of a fibre bundle with ${\bf R}^p-\{0\}$ as base.  
 To show that the period space is homotopic to $S^{p-1}$, it is sufficient
 to show that every fibre is contractible. This is an easy exercise that
we leave
 to the readers.

 As the period space is homotopic to $S^{p-1}$, there must be some $(p-1)-$cycles which 
generate the $(p-1)$st  homology. 
It is useful to exhibit  such cycles; for this purpose, let us present another
picture of the period space.

The period space 
projects down to $V$. The preimage of an element $x\in V$ is
the set of positive $p-$planes which do not lie in $x^{\perp}$, where 
 $x^{\perp}$ denotes the orthogonal complement of $x$ in $V$. 
Now decompose the vector
 space $V$  into three parts, $V=V_+\cup V_-\cup V_0$, with
$$V_+=\{x\in V|x\cdot x>0\},\quad V_-=\{x\in V|x\cdot x<0\},\quad
V_0=\{x\in V|x\cdot x=0\}.$$
$V_0$, usually called the light cone,
 is diffeomorphic to ${\bf R}\times S^{p-1}\times S^{q-1}$.   $V_+$ is diffeomorphic to 
${\bf R}^+\times S^{p-1}\times B^q$ and $V_-$ is diffeomorphic to 
${\bf R}^+\times B^p\times S^{q-1}$ and they intersect along $V_0$.
If $x\in V_+\cup V_0$, no positive $p-$plane can be
perpendicular to it. 
 Therefore, over $V_+\cup V_0$, the period space has a 
fibre bundle structure with contractible fibres, $G_p^+(V)$. 
If $x\in V_-$, the $p-$planes perpendicular
 to it are the positive $p-$planes which do not lie in $x^{\perp}$ 
and is itself a sub-self-dual Grassmannian.  We denote it by
$G_p^+(x^{\perp})$.  $G_p^+(x^{\perp})$ is a submanifold in $ G_p^+(V)$ 
 of codimension $p$. In fact, following the same line of the previous 
argument, one can easily see that $G_p^+(V)$
 is diffeomorphic to an affine space bundle over $G_p^+(x^{\perp})$ of fibre
dimension $p$.  From here it follows again that the homotopy type of 
$G_p^+(V)-G_p^+(x^{\perp})$ is just $S^{p-1}$, a sphere of dimension $p-1$. 
Let us denote the complement of $G_p^+(x^{\perp})$ in $G_p^+(V)$
by $\hat G_p^+(V,x)$.
 Any cycle which generates the $(p-1)$st homology of 
$\hat G_p^+(V,x)$ over $x\in V_-$ is the homology generator of the period 
space.
On the other hand, the 
$S^{p-1}$ factor in $V_+$ would also 
be the generator of the $(p-1)$st homology of the 
period space. This simple observation, though elementary, will 
play a very 
crucial role in the later development of the whole theory.

 Let us summarize the picture when $B$ is a point. The
period space of  $V=H^2(M;{\bf R})$ is homotopic to the sphere
$S^{p-1}$, where $p$ is equal to $b^+$. The number of ${\cal L}$
chambers is greater than one only if $b^+=1$.  
And when $b^+=1$,  the corresponding 
 period space is homotopic to $S^0$, which
consists of two points. In fact this is why the usual chamber structure
  of the ordinary Seiberg-Witten theory has a ${\bf Z}_2$ grading. 

For general $B$,
 the ${\bf Z}_2$ grading is not enough to parametrize the chambers.
 In fact we will see shortly that
it is very typical that the chambers are at least graded by $\bf Z$.
 Namely, there are an infinite number of chambers for even a single
Spin$^c$ structure. This would be a general characteristic difference
from the ordinary Seiberg-Witten theory.

 Suppose we are given a fibre bundle ${\cal X}$ over $B$ whose
 fibres are diffeomorphic to
$M$. The cohomology of the fibres $H^2({\cal X}_b; {\bf R})$  forms a flat 
vector bundle ${\cal V}$ over the base $B$. 
We will define the period bundle via ${\cal V}$,  
which generalizes the period space
and plays the key role in understanding the chamber structure.

\proclaim Definition 3.4. The period bundle ${\cal P}\longrightarrow B$ is defined to be the
fibre bundle canonically constructed by the flat bundle 
${\cal V}=H^2({\cal X}_b;{\bf R})
\longrightarrow B$ through the period construction. \par

The period bundle is a subbundle of
$\tilde{\cal P}=G_p({\cal V})\times_B {\cal V}$ 
and has a canonical surjective bundle morphism to ${\cal V}$.

From Proposition 3.2, we immediately have 

\proclaim Proposition 3.5. The homotopy type of the period bundle:  The period 
bundle ${\cal P}$ is homotopic equivalent to an $S^{p-1}$ bundle over 
$B$ through fibrewise homotopy. \par

Let us end this subsection with a remark concerning a result of Kronheimer.

\noindent{\bf Remark 3.6}. Given a  four-manifold $M$ with a symplectic form
$\omega_0$,
consider the space 
$$\Lambda_0=\{\omega\in\Omega^2(M)|
\omega\hbox{ is symplectic and cohomologous to }\omega_0\}.$$
 For each Spin$^c$ structure determined by an element 
$e\in H^2(M;{\bf Z})$ whose index 
$$d=K_{\omega_0}\cdot e-e\cdot e$$
lies in the range $0< d< b^+-1$ and which satisfies 
$\omega_0\cdot e<0$, Kronheimer [K] defined a homomorphism
$$Q:H_{d-1}(\Lambda_0;{\bf Z})\longrightarrow {\bf Z}$$
via Seiberg-Witten equations for the Spin$^c$ structure.

From Proposition 3.3, 
we have
$$H_i(\Gamma-{\cal WA}; {\bf Z})=H_i({\cal P}; {\bf Z})=H_i(S^{b^+-1};{\bf Z}).$$
Tracing Kronheimer's definition of $Q$, 
it is not hard to see that $Q$ can be extended to
 any Spin$^c$ structure determined by an element 
$e\in H^2(M;{\bf Z})$ whose index lies in the range $d>b^+$ and which satisfies 
$\omega_0\cdot e<0$.

This extension is particularly interesting in the case $b^+$ is equal to
one. 

\medskip
\noindent{\bf \S3.2. The chamber structure}
\medskip
In this subsection, we will describe the chamber structure in 
cohomological terms. 

Fix a Spin$^c$ structure ${\cal L}$ with $L$ as the associated
complex line bundle. The map $pr_{\cal L}$ maps $\Gamma-{\cal WA}({\cal L})$
onto sections of ${\cal P}$. Furthermore,
two pairs in $\Gamma-{\cal WA}({\cal L})$ lie in the same ${\cal L}$ chamber
if and only if their images under $pr_{\cal L}$ are homotopic as sections of 
the period bundle ${\cal P}$. Therefore we have  

\proclaim Proposition 3.6. Chamber structure: The chambers of a single 
Spin$^c$ structure are classified by fibrewise homotopy classes of sections
 into ${\cal P}$, denoted by $[B,{\cal P}]_f$. 


Unlike the case  $b^+=1$,
 the chamber structure  can be extremely complicated. 
For example
 consider the case $B=S^a$ and the fibre bundle ${\cal X}$ a 
trivial product $M\times B$ with $b^+(M)=b$ and $ a>b$. Then the set of
fibrewise
homotopy classes $[B,{\cal P}]_f$ is given by 
the set of free homotopy classes
 from $S^a$ to $S^b$, $[S^a,S^b]$,
a very complicated object in homotopy theory.
 In general if we consider a product fibre bundle ${\cal X}$,
 the set $[B,{\cal P}]_f$ will be given by $[B,S^{b^+(M)-1}]$, the 
cohomotopy set of $B$,     
 and it does not
 have a group structure.
We do not plan to give a complete description
of these chambers. Instead, we would only like to discuss how to measure the
difference of two chambers in cohomological terms. This piece of information
 will play a crucial role in the derivation of the wall crossing formula.

Whether two sections of  $\cal P$ are fibrewise homotopic
 to each other is naturally an extension 
problem: whether the  map
$f:B\times\{0,1\}\longrightarrow {\cal P}$ given by
the two sections can be extended to a map $F:B\times [0,1]\longrightarrow {\cal P}$.
This can be achieved via  obstruction theory.
 The answer to this question is that there will be
a sequence of cohomology classes associated to it which measure
 the obstruction of 
extension order by order. More precisely, we have ([S])

\proclaim Proposition 3.7. The obstruction to extending the map $f$ to $F$ is
measured by a sequence of elements which live in 
$H^r(B;\pi_r^{local}({\cal P}))$ 
for all $r$. Here $\pi_r^{local}({\cal P})$ is the local system of $r-th$ 
homotopy groups of fibres.
 In particular, if all these obstruction classes vanish, then the
 extension $F$ exists.\par

 As we know that the period space is homotopic to $S^{b^+-1}$, the 
obstruction classes mentioned in the proposition actually lie in 
$H^r(B;\pi_r^{local}(S^{b^+-1}))$ for all $r$.   
 In general the coefficient groups are not the integers. However there is 
one element which is particularly important, the principal obstruction.
 It is the first obstruction for the extension, living in $H^{b^+-1}(B;{\bf Z})$.
 Let us denote it by $Obs(s_1,s_2)$ for $[s_1],[s_2]\in [B,{\cal P}]_f$.
Even though the homotopy classes themselves
 do not have any additive property, the obstruction class discussed here 
does have the following remarkable property:
$$Obs(s_1,s_2)+Obs(s_2,s_3)=Obs(s_1,s_3).$$

It is this cohomology class $Obs$ which will enter the formulation of
the  wall crossing formula. The subsequent discussions will
illustrate this point.
Suppose we are given a generic one-parameter family of 
fibre metrics and two-forms, parametrized by $I=[0,1]$, over ${\cal X}/B$,
such that, at $0$ and $1$, no reducible solutions occur. 
Let ${\cal S}$ denote the subset in $B\times I$ where reducible solutions
occur.
 The singular locus
 $\cal S$ is usually complicated and depends on the 
one-parameter family of
metrics and the given two-forms. 
But what is important to us is the homology class of the wall crossing
locus ${\cal S}$ and this is described by 
\proclaim Lemma 3.8. 
The singular locus $\cal S$ represents a homology 
class of $B\times I$ in degree $b^+-dimB-1$.
Under the natural identification between 
the homology of $B$ and the homology of $B\times I$, this class
 is Poincar\'e dual to 
the principal obstruction class $Obs$.\par   

\noindent {\it Proof}. 
Under the map $pr_{\cal L}$, the one-parameter
family of fibre metrics and two-forms  gives a section of
${\cal P}^+$ (precisely, the pull back of 
${\cal P}^+$ from $B$ to $B\times I$) 
over $B\times I$.  ${\cal S}$ is just the 
singular locus in $B\times I$ where the section lies outside of ${\cal P}$
(also the pull-back).
 By transversality we can always perturb the one-parameter
family such
that the family is transverse, so $\cal S$ is the 
transversal intersection between the nonzero section and the zero section.
such that $\cal S$ 
 is embedded into $B\times I$ as a 
submanifold.

As the original extension problem is homotopic in nature,
we can assume our period bundle is an ${\bf R}^{b^+}-\{0\}$ bundle over 
$ B$. 
 To relate the obstruction class  
with $\cal S$, we need a CW complex structure on 
 the manifold with boundary
$B\times I$ which is transversal to $\cal S$. 
One way to obtain such a CW complex structure is
to  triangulate $B\times I$ such that 
the triangulation restricted to $\cal S$ gives a
 triangulation of $\cal S$ as a smooth submanifold and 
take the dual CW complex associated
 to this particular triangulation.  
Assume a CW complex structure transversal to $\cal S$ is given, we 
want to use it to construct the obstruction class. In
 this case the map has been defined
 outside  $\cal S$. Thus, over the cells which do not intersect $\cal S$ at all,
 there is no obstruction to extending the map from low strata. If we extend
 the map from the zero stratum to the higher ones, the first time the 
obstruction appears is in dimension $b^+$.  
The map still extends automatically across
 any $b^+-$dimensional cell which is disjoint from $\cal S$ 
and those cells are assigned to the 
zero element in $\pi_{b^+-1}(S^{b^+-1})$. On the other hand, any 
$b^+-$cell which intersects $\cal S$ is assigned to $\pm1$
 in $\pi_{b^+-1}(S^{b^+-1})$ by the transversality assumption.
Whether it is
  assigned to $1$ or $-1$ is  according to the orientation. At the end, we
 get a $b^+-$cochain of the CW complex.   It is a cocycle on $B\times I$
representing the principal obstruction class. On the
other hand, from the Poincar\'e duality,
it  is exactly the Poincar\'e dual to $\cal S$. This ends the
proof of the lemma. 

   As just mentioned, this elementary result 
will be used in the derivation of the wall crossing
 formula later. The readers should keep in mind that it is the cohomology
 class associated to $\cal S$, not the manifold $\cal S$ itself, which is
relevant to our later discussion.

\noindent{\bf  Example 3.9} (winding family). Let $M$ be a manifold
with $b^+>1$ and ${\cal X}$ be a
trivial bundle with base $S^{b^+-1}$. Then there is a
canonical chamber characterized by containing the constant parameters,
i.e. parameters of constant metrics and constant self-dual
two-forms. With the canonical chamber understood, the 
chambers are naturally identified with the integer group ${\bf Z}$.
A family of closed two-forms with positive square is called
a winding  family of two-forms if they represent a 
generator of the $(b^+-1)$st homology of $V^+$ introduced in the 
last subsection, and it is called a symplectic winding 
family if all the two-forms are symplectic.
It is not hard to see that the chamber containing a parameter of constant metric and 
winding two-forms corresponds to $\pm 1$ under the aforementioned 
identification. 

Interesting symplectic winding families include hyperk\"ahler families of 
K3, the four-torus and some $S^1$ families of primary Kodaira surfaces
studied in [Ge].

\bigskip
\noindent{\bf \S4. The wall crossing formula}
\medskip

 In this section we would like to discuss the wall crossing formula for 
the family version of Seiberg-Witten invariants. As before let 
${\cal X}\longrightarrow B$ be a fibre bundle of four-manifolds. 

In the first subsection, 
 we would like to discuss the phenomena of wall crossing
 in some generality. In the second subsection, 
 we would like to specialize to the case of 
$b^+-1=dim B$ where we can derive 
a universal wall crossing formula. 
We call this special case the ``critical'' case.
In the third subsection, we discuss the Fulton-MacPherson spaces.

\medskip
\noindent{\bf \S4.1. The Kuranishi models}
\medskip
 In the case $b^+-1> dim B$, we can always move the metrics and two
forms perturbation in a one-parameter family and show by the standard 
bordism argument that the invariants are actually independent of  the
 metrics and two-form perturbations. However ,
 in the $b^+-1\leq dim B$ case, the general position argument
 cannot provide us a completely smooth cobordism between these two
 moduli spaces. Typically, the
parameters (including the metrics and two-forms) will hit the singular values
 somewhere 
and reducible solutions develop in the bordism. In this case, 
the bordism is no longer a smooth bordism between these 
two different moduli 
spaces.  To understand how  the invariance property of the ``invariants''
 fails, we need to understand the behavior of the moduli space near the
reducible solutions. The following is the fundamental tool to describe the
 neighborhood of the moduli space near a solution $(A_b,\psi_b)$. 

\proclaim Proposition 4.1 (Family Kuranishi Model). Let $(A_b,\psi_b)$ be a 
solution
of the family Seiberg-Witten equations at $b$; then a neighborhood of the moduli 
space near this point $(A_b,\psi_b)$ can be described by $F^{-1}(0)/S^1$, where
$$F:T_b B\times H^1({\cal X}_b;{\bf R})\times U_b
\longrightarrow \hbox{Ker}(d^{\ast})\times \hbox{Coker}(D_{A_b})$$
 is an  $S^1$ equivariant map and $U_b$ is a small ball in 
$ Ker(D_{A_b})$.

\par

 The proof of this proposition is standard (see [DK]).
 When the spinor $\psi_b$ is not identically zero on the manifold ${\cal X}_b$,
 the $S^1$ action is free and the Kuranishi model does not tell us
too much other than that it is
 a smooth point in the moduli space. However when the
spinor $\psi_b$ is identically zero,  $S^1$ acts trivially on the point
 $(A_b,0)$ and the Kuranishi model describes  the local 
singular behavior near this point.   To get useful global topological 
description for the purpose of calculation, 
we actually need a global version of the
 Kuranishi model 
which generalizes the picture of the $b^+=1 $ case in [LL1]. Let us
 formulate it as the Global Kuranishi model for the bordism.

 Suppose there is a one-parameter family of deformation of the parameters
 $t, 0\leq t\leq 1$. Let ${\cal B}$ denote the bordism connecting the
 two different moduli spaces at $t=0$ and $t=1$. Suppose that somewhere
 between $0$ and $1$ the bordism hits the wall, then we have

\proclaim Proposition 4.2 (Global Kuranishi Model). Let ${\cal S}$ denote the
reducible solutions in the bordism ${\cal B}$ connecting two different smooth 
moduli spaces, and  denote the projection map of ${\cal S}$ to $I\times B$
by $\pi$ and the image by $S_0$. Suppose that ${\cal S}$ and ${\cal S}_0$
 are both smooth, 
then the local neighborhood of ${\cal S}$ in ${\cal B}$ 
is described by the following Kuranishi model $F^{-1}(0)/S^1$, where
$$F: T{\cal S}\times \pi^{\ast}N_{{I \times B}/{\cal S}_0}\times
U \longrightarrow 
\hbox{Ker}(d^{\ast})\times \hbox{Coker}(D_A)
$$  
 is again an $S^1$ equivariant map, $N_{{I \times B}/{\cal S}_0}$
is the normal bundle of ${\cal S}_0$ inside $I\times B$,
$U$ is a small ball in $ \hbox{Ker}(D_A)$ and
$A$ varies over the reducible connections in $\cal S$. \par

However, this global model still has limited use.
 Usually we do not know exactly the diffeomorphism type of the
 singular set ${\cal S}$,
nor do we know that ${\cal S}$ and ${\cal S}_0$ are smooth.
 To get a more useful model, 
 it must be replaced by a fattened version.

Let  ${\cal T}^{b_1}_B$ be the $T^{b_1}$ bundle over $B$ which parametrizes
 the reducible connections in the family ${\cal X}\longrightarrow B$.
It is well known that the dimensions of the kernel spaces and
cokernel spaces may vary on ${\cal T}^{b_1}_B$ and 
$ \hbox{Ker}_B(D_A)$ and $ \hbox{Coker}_B(D_A )$ are in general not
honest complex vector bundles over ${\cal T}^{b_1}_B$. To avoid this difficulty, we can use the technique 
 in [LL1] to perturb the Kuranishi model a little bit (see also
the technique developed in [R]). In the present case,
the technique can be applied as
the base ${\cal T}^{b_1}_B$ is still compact. And after this slight
perturbation, we can always assume $\hbox{Ker}_B(D_A)$ and $\hbox{Coker}_B(D_A)$ 
are honest vector
 bundles. As in [LL1], we have

\proclaim Proposition 4.3 (Modified Global Kuranishi Model).
There exists an $S^1$ 
equivariant, fibre preserving map
$F$ from a disc subbundle in the complex
vector bundle $ \hbox{Ker}_B(D_A)$ to the complex vector bundle
$ \hbox{Coker}_B(D_A)$ whose zero set contains the neighborhood of
${\cal S}$ in ${\cal B}$ as an embedded subspace.
\par

When $B$ is 
 a single point, the picture reduces to the previous one
described in [LL1]. From now on, let us drop the subscript $A$ and write $D_A$ simply as $D$.

 By the modified global 
modified Kuranishi model, the link space of the singularities 
$Link({\cal S})$ in ${\cal B}$ can be
 always embedded inside the projective space bundle formed by 
$\hbox{Ker}_B(D)$. We do not really care about its topology in detail, 
in fact  we  do not even assume that ${\cal S}$ and ${\cal S}_0$
are smooth. What we only need 
to know is 
 the homology  class represented by $Link({\cal S})$. It is Poincar\'e dual to a certain obstruction class
   over $P(\hbox{Ker}_B(D))$.
 However there is some slight difference from the usual $b^+=1$ wall crossing
 formula.
 Namely, there are two sources of the obstruction classes. One of them is
comes from the obstruction bundle which is described in the
following definition. 

\proclaim Definition 4.4 [LL1]. Obstruction Bundle ${\cal O}bs$.
Let us denote the projection map from $P(\hbox{Ker}_B(D))$ to ${\cal T}^{b_1}_B$ by $\pi_1$. The obstruction bundle ${\cal O}bs$
is given by
$${\cal O}bs=\pi_1^{\ast}\hbox{Coker}(D)\otimes H,$$
where $H$ is the hyperplane line bundle canonically associated with $P(\hbox{Ker}_B(D))$.
\par
 
 As argued in Lemma 2.2 in [LL1], $F^{-1}(0)$ can be viewed as
the zero set of a section of ${\cal O}bs$.

 Another source of the obstruction class comes from 
the period bundle and
 has been explained in detail in section 3.  
Suppose that the period maps at $0 $ and $1 $  both map into ${\cal P}$;
 then 
 in general there would be some obstruction to extending
 these into a well 
defined map from $B\times I$ to the period space. The so-called primary
obstruction  is measured by a cohomology class $Obs$ which lives in
$H^{b^+-1}(B;\pi_{b^+-1}{\cal P}_b)$. The coefficient is the local system
 formed by the $(b^+-1)$st homotopy group of the fibre of the period bundle.
 As discussed in  section 3,  the period space is homotopic 
to $S^{b^+-1}$. Thus the obstruction class $Obs$ lies in the cohomology
with integer coefficients.
 Then $Obs$ is Poincar\'e dual to the base locus 
in $B\times I$ where 
 the bordism hits the wall.

Let us state the result for a general invariant $SW({\cal X}/B, 
{\cal L}, \Theta)$. From now on, we will drop ${\cal X}/B$
in the notation and simply write the invariant $SW({\cal X}/B, 
{\cal L}, \Theta)$ as $SW({\cal L}, \Theta)$.
 
\proclaim Proposition 4.5.  The wall crossing number is calculated by
$$\pm WCN({\cal L},\Theta)=\int_{P(\hbox{Ker}_B(D))} 
\Theta\cdot H^{dim(\hbox{Coker}_B(D))}\cdot Obs, \eqno (4.1)$$
where the class $Obs$ is pulled back from the base $B$ to $P(\hbox{Ker}_B(D))$ and is 
viewed as a cohomology class living in the projective bundle.\par

{\it Proof}. First,
$$\pm WCN({\cal L},\Theta)= {\int_{Link({\cal S})} \Theta}.$$
Following [LL1], we  reduce
 the  calculation of the wall crossing number to the calculation on ${\cal T}^{b_1}_B$, 
with the insertion of the Euler class of the obstruction bundle 
${\cal O}bs$ and the
homotopy theoretical obstruction class $Obs$. Since we can assume, 
as explained in [LL1], that
$\hbox{Coker}_B(D)$ is a trivial bundle, we have
$$\eqalign {\int_{Link({\cal S})} \Theta=&\int_{P(\hbox{Ker}_B(D))}
\Theta\cdot Obs\cdot c_{dim(\hbox{Coker}_B(
D))}({\cal O}bs)\cr
=&\int_{P(\hbox{Ker}_B(D))} \Theta\cdot H^{dim(\hbox{Coker}_B(D))}\cdot Obs.\cr}$$

For the pure invariant, (4.1) takes the following simple form:
$$\int_{P(\hbox{Ker}_B(D))}  H^{dim B+b_1+ dim(\hbox{Ker}_B(D))-b^+ +1}
\cdot Obs.$$

Notice that whenever the base dimension is lower than $b^+-1$ then the
obstruction class $Obs$ vanishes automatically and the wall crossing number 
is zero. This statement coincides with the previous 
observation that the wall crossing phenomenon does not exist whenever
 $dim B<b^+-1$. 
Another simple fact is that 
if we consider a weight $m$ mixed Seiberg-Witten invariant and 
$m>dim B-b^++1$, then the wall crossing number is always zero.
This is because if the degree of the 
cohomology classes from the
base will exceed the base dimension $dim B$ and the invariants
 automatically vanish. 

In general, to calculate the integral (4.1) we need information
about the index bundle $\hbox{IND}_B(D)=\hbox{Ker}_B(D)-\hbox{Coker}_B(D)\in 
K({\cal T}^{b_1}_B)$; to be precise, we want to calculate
its Chern character  
$Ch(\hbox{IND}_B(D))\in H^{\ast}( {\cal T}^{b_1}_B;{\bf Q})$.
 It is calculated by
 the Family Index theorem as a push-forward from ${\cal X}\times_B {\cal T}^{b_1}_B$ 
to ${\cal T}^{b_1}_B$,
$$Ch(IND_B(D))=\int_{{\cal X}\times_B {\cal T}^{b_1}_B/{\cal T}^{b_1}_B}
{\cal A}_{{\cal X}\times_B {\cal T}^{b_1}_B/{\cal T}^{b_1}_B}
\cdot Ch({\cal L}\otimes \Omega),$$
where ${\cal A_{M/N}}$ denotes the relative $\hat A$ genus of the fibre 
bundle and $\Omega$ is the tautological line bundle on ${\cal X}\times_B
{\cal T}^{b_1}_B$. When we restrict the line bundle $\Omega$ to each 
fibre ${\cal X}_b\times {\cal T}^{b_1}_B$, its first Chern class is 
given by the standard formula $c_1(\Omega_b)=\sum x_i\cdot y_i$, where
 the $x_i,y_i$ are some dual bases of $H_1({\cal X}_b;{\bf Z}),
 H^1({\cal X}_b;{\bf Z})$,
respectively.

\bigskip
\noindent{\bf \S4.2. The critical case}
\medskip
In this subsection, we will 
concentrate on the critical case where $b^+-1=dim B$. 
In this case one has an explicit expression for the
wall crossing number because of the following lemma.

\proclaim Lemma 4.6. Suppose $b^+-1$ is equal to the dimension of $B$.
 Then the obstruction class $Obs$ is a top dimensional class in 
$H^{dim B}(B;{\bf Z})$. 

$Obs$ will produce a number which can be interpreted as the number
 of times that the parametrized moduli space  
intersects the wall when the intersections are transverse.
And when the intersections are transverse, we can think of ${\cal S}_0$ as simply a finite number of points.

Let us restrict to the case that the two chambers are consecutive. 
In this case, $Obs$ reduces ${\cal S}_0$ to a single point, 
and the link of
the reducible solutions is simply
 the traditional picture $P(\hbox{Ker} (D))$ of the projective
space bundle over $T^{b_1}$ in [LL1].

Notice that the wall crossing numbers between any two consecutive chambers
are the same. This simple fact will be crucial to deducing a vanishing 
result in the next subsection. 

Let us first derive the wall crossing formula for the pure invariant.
 By the above discussion, (4.1) is simply
 $$\eqalign{\int_{P(\hbox{Ker}(D))} H^{top}\cr} \eqno (4.2)$$
What we have to 
calculate is the top power of the hyperplane class on the projective
 bundle of the index bundle 
(we know we can choose the representative of the $K-$theory class
so that $\hbox{Ker}(D)-\hbox{Coker}(D)$ is an honest vector bundle).
In [LL1], we did some lengthy calculation  in the $b^+=1$ case.
 Actually the calculation can be highly simplified by noticing that
we calculated the $n-$th Segre class of the Index bundle
with $n=b_1/2$. Denote the $n-$th Segre class by $s_n$. 
By the same token,  what we want to calculate here is simply $s_n(\hbox{IND} D)$
with $n=dim B+b_1/2$. 
 The Segre classes and the Chern classes can be related to each other 
 in a simple way.
Namely if one defines $S(t)= 
\sum_i s_i t^i$ and $C(t)=
\sum_i c_i t^i$( where $c_i$ and $s_i$ denote the $i-$th chern class and
 Segre class respectively), then (see [F])
$$C(t) S(t)=1,$$
i.e, $S(t)$ is the inverse power series of $C(t)$! This observation leads to a
 simpler computation.

  The index calculation is performed
 on $M\times T^{b_1}$.
It follows from the family index theorem that the grade $2i$ term of the 
expression 
$$\int_{M}A(M)ch(L)ch(\Omega)$$
(where $\Omega
\in H^2(M\times T^{b_1};{\bf Z})$ 
denotes the first Chern class of the Poincar\'e
 line bundle) is the $i-$th Chern character, denoted by $ ch_i$. 
 The first term of the index calculation gives the formal dimension.

The following lemma calculates the Chern character.

\proclaim Lemma 4.7. $ch_i=0$, for all $i>2$ and
$$\eqalign{ch_1=&c_1({\cal L})\cdot \Omega^2/2!
[M]\cr ch_2=& 
\Omega^4/4![M].\cr}$$

\par

\noindent{\it Proof}. 
 Expanding $ch(\Omega)$ one gets various powers of $\Omega$.
From the fact that $\Omega$ is of bi-degree
(1,1) on $M\times T^{b_1}$, it cannot be raised to  more than fourth
power or it will automatically vanish
along $M$. 
 
 Notice that this is slightly different from the usual $b^+=1$ case
 where we automatically have the vanishing of $ch_i$ with $i\geq 2$
 because of the light cone lemma.

Let us temporarily forget about the detailed expression of $ch_1$ and $ ch_2$
 in terms of $\Omega$ and ${\cal L}$, etc.
It makes things simpler to keep the calculation symbolic.

The next step is to prove the following recursion relation between 
 the various $c_i$.

\proclaim Lemma 4.8. The Chern classes of the index bundle satisfy the
 recursion formula

$$c_i={1\over i} (c_{i-1}ch_1-2 c_{i-2}ch_2).$$\par

\noindent{\it Proof}. Let us look at the defining equation of the
cohomology of the projective bundle (the splitting space).
$$ x^m-c_1 x^{m-1}+c_2 x^{m-2}-...+(-1)^m c_m=0,$$
 where $m$ is the complex dimension of the bundle.
Notice that the variable $x$ can represent any of the Chern roots of the bundle
$\lambda_i$, $i=1,\cdots, m$. Thus we can
 plug all the Chern roots in to get $m$ equations. Summing these $m$ equations
 we get 
$$(\sum_i \lambda_i^m)-c_1(\sum_i \lambda_i^{m-1})+...+(-1)^m mc_m=0.$$
Using the fact $ch_i=0, i>2$ and $ch_1=\sum_i \lambda_i$, $2ch_2=\sum_i 
\lambda_i^2$
we see that 
$$mc_m=c_{m-1}ch_1-2c_{m-2}ch_2.$$
 This proves the recursion
 formula for $i=m$.
 For general $i$, one divides the equation by $x^k$ and  get
$$x^{m-k}-c_1 x^{m-1-k}+...(-1)^mc_m x^{-k}=0.$$
 Plugging in the Chern roots and summing over $i$ we get an equation which can 
be reduced to the following form after using $ch_i=0, i>2$.

$$\sum_{l,l<k}(-1)^{l-k}(c_{m-l}(\sum_i \lambda^{l-k}_i))+mc_{m-k}-c_{m-k-1}ch_1
+2c_{m-k-2}ch_2=0. \eqno (4.3)$$
To show that the recursion formula is true, one only needs to show that 
 the first sum over $l$ is equal to $-kc_{m-k}$.
 
To show this let us play a trick. Denote 
$\sum_{p>0}\sum_i (\lambda_i t)^{-p}$ by $F(t)$, a Laurent series in
 $t$. Then the sum in (4.3)
 can be 
expressed as the coefficient of the $ (m-k)$th power  of $t$ of $C(t)\cdot F(t)$.
On the other hand, if we commute the sum over $p$ and $i$ in the definition of
$F(t)$ and extend $p$ to $\geq 0$, then we find that $F(t)$ can be rewritten as
$$\eqalign{\sum_i(\sum_{p\geq 0}(\lambda_i t)^{-p}-1)
=&\sum_i((1+(\lambda_i t)^{-1})^{-1}-1)\cr
=&-\sum_i((1+\lambda_i t)^{-1}).\cr}$$

Thus $C(t)F(t)$ can be simplified as $$-\prod_i (1+\lambda_i t)  (\sum_i
(1+\lambda_it)^{-1})=-\sum_i(\prod_{j\not= i}(1+\lambda_jt)).$$
 The final expression is related to taking the derivative of a product of
linear factors. In fact let us call $G(s)=s^m C(s^{-1})=\sum_r c_r s^{m-r}$.
 Then $-s^{1-m}(d/ds)(G(s))|_{s=t^{-1}}$ is exactly equal to this 
expression. Thus it is easy to see that 
the $(m- k)-$th order term in $t$ is exactly 
$-kc_{m-k}$. This proves the lemma.

Knowing the recursion relation, it follows that one can derive an
O.D.E.
 for $C(t)$.
More precisely, view $C(t)$ as 
an element in $R[[\lambda_1,\lambda_2,
\cdots,\lambda_m,t]]$, then we have 

\proclaim Lemma 4.9. The series $C(t)$ formally satisfies the following 
differential equation,
$${d\over dt}C(t)=(ch_1-2ch_2t)C(t).$$

On the other hand we know that $C(0)=1$. With this 
 initial condition on
 the differential equation,  we immediately get 
$$C(t)=Exp(ch_1 t- ch_2 t^2).$$
Therefore, by $C(t)S(t)=1$, 
$$S(t)=Exp(-ch_1 t+ch_2 t^2). \eqno (4.4)$$

To get the $n-$th order term of $S(t)$ we collect the terms on the
right hand side  of (4.4) which are of the type $t^i t^{2j}$ such that $i+2j=n$.
Each term of this type contributes  $(-1)^i ch_1^i ch_2^j/i!j!$,
 thus we get our
final formula.

\proclaim Theorem 4.10. The wall crossing number of the pure invariant 
in  the critical case is
expressed in term of the Chern characters as follows,

$$\pm WCN({\cal L}, H)=\sum_{i+2j=b_1/2}(-1)^i ch_1^i \cdot ch_2^j/
i!j![T^{b_1}].
\eqno (4.5)$$
\par

This theorem applies to section invariants as well
and we state it as a proposition.

\proclaim Proposition 4.11. The wall crossing number of
section invariants is also given by (4.5).

{\it Proof}.
Though different sections $s$
give rise to different cohomology classes $e(s)$ on the moduli space
which is a fibre space over $B$,
they all determine the same cohomology class on $P(Ker(D))$, and this class
is nothing but the hyperplane class $H$. 
 Hence they all have the same wall crossing number as $WCN({\cal L},H)$.

If one uses the fact that $ch_2=0$ for $b^+=1$ four-manifolds, one
 immediately recovers the wall crossing formula in [LL1].
 The explicit formula of $ch_i$ depends on some details of the  structure of the
cohomology ring  of $M$. However, when 
 $M$ has $b_1=0$, the wall
 crossing number of the pure invariant is simply given by $\pm 1$.
A particular interesting example is
an $S^2$ winding family of K3 surfaces (recall that  a K3 surface has
$b^+=3$).

For any manifold with $b_1\not=0$, we would like to describe the form 
of its wall crossing formula in terms of its cohomology ring pairing.
 First we know that $ch_1$ is calculated by $c_1({\cal L})\cdot \Omega^2/2!
[M]\in H^2(T^{b_1})$ and $ch_2$ is calculated by 
$\Omega^4/4![M] \in H^4(T^{b_1})$.
 Let us define ${b_1\choose 2}$  linear functionals $ q_{ij}$ on 
$H^2(M)$ by 
 cupping $c_1({\cal L})$ with $x_i\cdot x_j$. Then $ch_1$ can be 
re-expressed as
 $\sum_{i<j}q_{ij}({\cal L})y_i\cdot y_j \in H^2(T^{b_1})$. 
 Likewise we can introduce ${b_1\choose 4}$ numbers, denoted by 
$\epsilon_{ijkl}$,
 by cupping $x_i,x_j,x_k,x_l$ together. Then $ch_2$ can be expressed 
as $\sum_{i<j<k<l}\epsilon_{ijkl}y_i\cdot y_j\cdot y_k\cdot y_l
\in H^4(T^{b_1})$.

Even though the formula looks complicated we are still able to 
 conclude that the
wall crossing formula is a polynomial function in $c_1({\cal L})$ of 
degree at most $b_1/2$. Whether this function is totally trivial or is 
honestly a polynomial of degree $b_1/2$ really depends on the ring 
structure of the cohomology of 
$M$. In case we know that the cup product of $H^1$ is 
nondegenerate, we always know that the wall crossing number behaves 
honestly as a  polynomial of degree $b_1/2$. This is the case for
K\"ahler surfaces where the hard Lefschetz theorem implies the 
non-degeneracy of the cup product. On the other hand, if the cup product 
on $H^1$ is totally null, it automatically implies the vanishing of 
 the wall crossing number of the section invariants of
all the Spin$^c$ structures.

Let us look at an important example,
any  $S^2$ family of the four-torus $T^4$ ($T^4$
also has $b^+=3$). Unlike the case of the K3 surface,
the wall crossing numbers of $T^4$ are no longer just $1$ or $-1$
and they 
 are expressed as a degree two polynomial in ${\cal L}$,
 as was expected from the K\"ahler property of $T^4$.
A simple calculation tells
 us that $ch_1=(c_1({\cal L})/2)^2$ and $ch_2=1$.
 Since $WCN({\cal L},H)=ch_1^2/2+ch_2$, we have

\proclaim Corollary 4.12. The wall crossing numbers of the pure invariant of a 
Spin$^c$ structure ${\cal L}$
 on an $S^2$ family of $T^4$ are given by 
 $$\pm  WCN({\cal L}, H)=\pm({1\over 8} c_1({\cal L})^2+1).$$
\par

Before going on to discuss the wall crossing formula for 
general Seiberg-Witten invariants 
in the critical case, we want to make a remark
concerning the noncritical case.

In the critical case, the formula is very neat 
because $ch_i=0,i>2$.
 However, even without this assumption,
  $c_i$ still satisfies the following recursion relation,
 $$c_i={1\over i}(c_{i-1}ch_1-2c_{i-2}ch_2+3!c_{i-3}ch_3+\cdots ),$$
 as can be derived along the lines of lemma 4.8.

 Similarly  $C(t)$ satisfies the following differential equation,
$${d\over dt}(C(t))=C(t) (ch_1-2!ch_2+3!ch_3-4!ch_4+\cdots).$$
 And the same initial condition leads to the solution
 $$C(t)=Exp(ch_1 t-ch_2 t^2+2!ch_3 t^3-\cdots)=Exp(\sum(-1)^{i-1} ch_i
t^i (i-1)!).$$
 The derivation is done along the same lines and we skip the proof.
 However the final expression for
 the Segre classes is much more lengthy than before.
  This formula has certain applications
 if we consider the family wall crossing formula of the noncritical cases.\par

Next let us discuss the invariants induced from $H_1(M;{\bf Z})$. 
In this situation, we must insert a certain power of the Euler class
 and a  certain 
power of the $\mu$ map images to the projective space bundle over $T^{b_1}$.
We  denote the class on $P(Ker (D))$ induced
 by $\mu(\zeta)$ also by $\mu(\zeta)$. It is well known that 
 this class further determines a class on the torus $T^{b_1}$. By further abusing notation,
 we denote this class by $\mu(\zeta)$ as well.
 As before, the effect of introducing the new classes merely reduces the
effective dimension of the torus $T^{b_1}$.

\proclaim Theorem 4.13. Given $\zeta_1, \cdots, \zeta_q
\in H_1$, the wall crossing number of $SW(\zeta_1, \cdots, 
\zeta_q)$
is given by
$$\pm WCN(\zeta_1,\cdots, \zeta_q)=\sum_{i+2j={b_1-q\over 2}}(-1)^i{i+j\choose i}ch_1^i \cdot ch_2^j/(i+j)!
\mu(\zeta_1)\cdots\mu(\zeta_q)[T^{b_1}]. \eqno (4.6)$$

 In general the wall crossing number will become a polynomial of $c_1({\cal L})$ 
of degree less than or equal to $b_1/2-q$.
  We would like to illustrate this by studying the primary 
Kodaira surfaces. 

The primary Kodaira surfaces are the complex
 surfaces (non-K\"ahler) which are quotients of ${\bf R}^4$ 
with $b^+=2$ and $b_1=3$. Let $(x,y,z,t)$ 
be the coordinates
 of ${\bf R}^4$, the discrete group is generated by $(x,y,z,t+1),(x,y,z+1,t),(x,y+1,z,t)$
and $(x+1,y,z+\lambda y,t)$ where $\lambda$ is a fixed positive number.
(We exclude the $\lambda=0$ case as it gives $T^4$). 
In symplectic geometry they are known to
 have been  constructed by Thurston as the first examples of 
symplectic manifolds which are not
K\"ahler.

If we project $M$ to $(x,t)$, $(y, t)$ or to $(x,y)$,
we obtain three descriptions of the primary Kodaira surfaces 
as a $T^2$ 
bundle over $T^2$ with $dx\wedge dt$, $dy\wedge dt$ and 
$dx\wedge dy$ as Poincar\'e duals to the fibres.

 
 

 It is easy to see that $H^1(M;{\bf R})$ is generated by the 
 three differential forms $dx,dy,dt$. On the other hand, $H^2(M;{\bf R})$ is 
generated by $dx\wedge dt,dy\wedge dt,dy\wedge dz,dx\wedge(dz-\lambda
xdy)$. From here it is easy to see that the subspace $dx\wedge dt,dy\wedge
 dt$ is the two-dimensional subspace of $H^2$ generated by $H^1$. 
 Comparing with $T^4$, one difference is that the ring structure
 on $H^1$ is not nondegenerate. There are two nice $S^1$ symplectic winding families,
$$\eqalign{\omega_0(\theta)=&cos\theta dx\wedge dt+sin\theta
dy\wedge dz\cr
\omega_1(\theta)=&cos\theta dy\wedge dt+sin \theta dx\wedge
(dz-\lambda xdy).\cr}$$
For generic $\theta$, the first fibration is symplectic with respect to
$\omega_0(\theta)$, and the second is symplectic with respect to $\omega_1(\theta)$.

The critical dimension is now $b^+-1=1$,
 thus it is natural to consider the family version of SW invariants with
 ${\cal X}=M\times S^1$.  
As the base is odd-dimensional, it is
 easy to deduce that the wall crossing numbers are all zero for
 the pure invariant! Thus we must
 consider invariants involving $H_1$. The same
observation also tells us that the insertion of an even number of classes in $H_1$
  always gives the answer zero. On the other hand, if we insert 
more than $b_1$  classes, it is also zero.  Therefore, we will 
only look at the cases $q=1$ and $q=3$.

 If we insert three classes into the definition of the invariants, it can be
 easily seen that the wall crossing numbers are given simply by
 $$\pm \mu(\zeta_1)\mu(\zeta_2)\mu(\zeta_3)[T^3].$$ 
 Whether it is zero or not depends only on whether $\zeta_i$ are linear
 independent in $H_1(M,{\bf R})$ or not. However the wall crossing formula
 can be applied only if the Spin$^c$ structure has a moduli space whose
 dimension is at least three.

 On the other hand, the $q=1$ case is less restrictive, the moduli spaces are
 only required to have dimensions bigger than zero. 

Any Spin$^c$ structure can be represented as a cohomology 
class $m_1 dx\wedge dt+m_2dy\wedge dt+m_3 dy\wedge dz+m_4 dx\wedge
(dz-\lambda xdy)$ with $m_i\equiv 0\pmod 2$ (We have used the fact that
 these differential forms are actually an integral basis of $H^2$).
 In the wall crossing formula, $ch_2\equiv 0$ as the cup product pairing
 on $H^1$ is degenerate.
 Thus the wall crossing formula is given by   
$$WCN({\cal L},\zeta)=-c_1({\cal L}/2)\cdot \Omega^2/2![M]\cdot \mu(\zeta)[T^3].$$
 
Denoting $\zeta^{\ast}$ by $n_1dt+n_2dy+n_3dx$, we find a numerical expression

$$WCN({\cal L},\zeta)={-(m_3\cdot n_2+m_4\cdot n_3)\over 2}.$$

Finally, to see the structure of the invariants, we can form their generating series. Let us define the generating series $WCN(\zeta)$ to be
$$\sum_{k_1,k_2,k_3,k_4}WCN(2(k_1dx\wedge dt+
k_2 dy\wedge dt+k_3dy\wedge dz+k_4dx\wedge( dz-\lambda dy)),\zeta)t_1^{k_1}
t_2^{k_2}t_3^{k_3}t_4^{k_4}.$$
With some simple reduction we get

$$WCN(\zeta)={-(n_2t_3+n_3t_4-(n_2+n_3)t_3t_4)\over 
(1-t_1)(1-t_2)(1-t_3)^2(1-t_4)^2}.$$

Similarly we define $WCN(\zeta_1,\zeta_2,\zeta_3)$ and it can be easily
computed  
$$\prod_{i=1}^{4} {1\over (1-t_i)}\phi(\zeta_1^{\ast}\wedge
\zeta_2^{\ast}\wedge\zeta_3^{\ast}),$$
 where $\phi:\wedge^3 H^1\longrightarrow {\bf R}$
 maps $dx\wedge dy\wedge dt$ to 1. 

It is easy to deduce that 
if $\cal L$ corresponds to a cohomology class not of
 the form $m_1 dx\wedge dt+m_2 dy\wedge dt$, then the Seiberg-Witten 
moduli space in the first winding chamber 
is nonempty  if $dim{\cal M}_{\cal L}\geq 1$.

\bigskip
\noindent{\bf \S4.2. Important examples of the non-critical case:
Fulton-MacPherson spaces}
\medskip

  In the last subsection, we  were able to derive the
 family wall crossing formula in the critical case. In that case, 
 the formula is universal in the sense that
it is  independent of the details regarding the
topology of the fibre bundle.  
When  $dimB$ is
 bigger than $b^+-1$, the picture is much more complicated.
In this subsection, we will  
thoroughly investigate an important class  of examples, the
Fulton-MacPherson spaces. 

The discussion of the Fulton-MacPherson spaces has a 
twofold purpose. On the one hand, it demonstrates the complexity
of the problem, namely that it is usually difficult to derive an explicit
 formula in the non-critical case.
In general,
 the answer and the detailed calculations depend on the topology of the
  given fibre bundle. As a result, we are not able to get a universal 
formula as we did before. 

On the other hand, it was suggested to the authors by
G. Tian  that nodal Gromov-Witten invariants should be
 differentiable invariants and that this might be proved using family Seiberg-Witten
 invariants over some natural families. We believe the
Fulton-MacPherson families are the right families.

The Fulton-MacPherson spaces will be used
 to build up a fibre bundle.
If we want to apply
 this scheme to 
K3, $T^4$ or the primary Kodaira surfaces, we will thicken the base by $S^2$ or
 $S^1$. In general, we will thicken the base by $S^{b^+-1}$.
This family is ultimately related to the counting of nodal (or 
higher singularities) pseudo-holomorphic curves. 
  In fact after certain transformation, the wall crossing formula in these
particular examples can be shown to be closely related to an enumeration
 problem in algebraic geometry.  This also strongly supports (as a nontrivial
 example) the general conjecture between the family Seiberg-Witten
 theory and the family 
Gromov theory. In fact, the result will be compared with the
 results from the blow-up formula by A.K. Liu in another up-coming paper.  
The fact that they coincide gives rise to 
 a nontrivial consistency check of the calculations. 

 First, let us introduce the concept of the Fulton-Macpherson spaces.
 Even though their construction is completely general, we will restrict 
our attention to almost complex four-manifolds. 
 Let $Y$ be an almost complex manifold. Let $Z$ be an almost complex 
submanifold inside $Y$ of codimension $m$. Then there is a well defined
 process of blowing up the submanifold $Z$ inside $Y$. Topologically, it
 is given by replacing an almost complex neighborhood of $Z$ by a 
 projective bundle over $Z$ such that the new manifold has a canonical
 almost complex structure. If $Y$ is an almost complex manifold with a
 symplectic structure, one can perform the symplectic blow-up by
 gluing the symplectic structure on the complement of the neighborhood and
 the symplectic structure on the projective space bundle by some cut-off
 process. For the details of the construction, consult Guillemin-
Sternberg and McDuff's papers [GS] and [M2].  
From now on  $M$ will be
assumed to be an almost complex four-manifold.

A  set of $n$ ordered  distinct points on $M$ is equivalent to a point in
the configuration space
$M^n(\triangle)=(M^n-diagonals)$, where
$M^n$ is the product of $n$ copies of $M$. 
The natural question to ask is whether
$M^n(\triangle)$ has 
 a smooth compactification which respects the action of the symmetric group 
$S_n$. It was Fulton and MacPherson ([FM]) who first 
explicitly constructed a good 
compactification  of the configuration space when $M$ is a 
 smooth complex manifold. Later a similar  construction
 was also used by Kontsevich and other people to study Chern-Simons theory.

 The Fulton-MacPherson space $M[n]$ is a nice compactification 
of 
$M^n(\triangle)$ such that the compactifying divisors have a certain moduli 
meaning. Moreover, the space $M[n]$ is smooth, $S_n$ equivariant
and maps
 surjectively to the space $M^n$. 
$M[n]$ can be described in several ways; we will construct $M[n]$ by a sequence
of $2^n-n-1$ blowups from $M^n$. For each of the $2^n-n-1$ subsets  $S$ of $\{1,\cdots, n\}$ with cardinality at least two, let $\Delta_S\subset M^n$
be the diagonal where the points $x_i$, for $i$ in $S$, coincide. 
We start with $M[2]$; it is just the blowup of $M\times M$ along the diagonal
which corresponds to the diagonal $\Delta_{\{1,2\}}$. Suppose
$M[n]$ is constructed with a map to $M^n$. We will construct
$M[n+1]$ by blowing up all diagonals $\Delta_S$ in $M^{n+1}$ 
(or more precisely, the preimage of the diagonals under the map
$M[n]\times M\longrightarrow M^{n+1}$) where $n+1$ is an element of $S$.
The blowups are done in $n$ steps. The first step is to
blow up the diagonal $\Delta_{\{1,\cdots, n\}}$, next we blow up
the proper transform of the $n-1$ diagonals $\Delta_S$ where $S$
has cardinality $n$. The $k-$th step is
to blow up the proper transform of the diagonals $\Delta_S$ where $S$ has cardinality $n+2-k$. 
Notice that after the $k-$th step, the proper transform of
the diagonals to be blown up in the $(k+1)$st step become
disjoint and the order of the blow-up of these
diagonals in this step is irrelevant. This is exactly
why $M[n]$ preserves the symmetry of the group $S_n$.
It is not hard to see that all the blowups are 
along submanifolds of complex codimension $2$ or $3$.

 Given this inductive description of $M[n]$, one may wonder
whether $M[n+1]\longrightarrow M[n]$ gives rise to a smooth fibre bundle
which parametrizes blowups.
Unfortunately some  fibres of this map are singular.
When a point lies in $M^n(\Delta)$, the fibres are all diffeomorphic
 to $M$ with
these $n$ points blown up. On the other hand, when two or more points 
collide,
 the fibres themselves will become unions of several normal crossing
 four-manifolds. One of them is $M$ with several points blown up,
the other
components are all rational (i.e. ${\bf CP}^2$ with fewer points blown up). They
intersect each other in ${\bf CP}^1$ and the total 
number of exceptional ${\bf CP}^1$s sums up to be
 $n$. As the current Seiberg-Witten theory is developed for  smooth
 fibre bundles, we restrict our discussion to this case (the Seiberg-Witten
 theory for singular spaces or singular fibrations is an interesting subject 
on its own).

However a simple modification will give us the desired bundle.
We start from the space $ M[n]\times M$ and denote it by ${\cal X}_0$.
 The composition of $M[n]\longrightarrow
M^n$ with the $i-$th ($1\leq i\leq n$) projection map
gives rise to $n$ graphs $V_i:M[n]\longrightarrow M[n]\times M$,
corresponding to the diagonals $\Delta_{\{1,\cdots, n+1\}},
\cdots, \Delta_{\{n,\cdots, n+1\}}$.  
Instead of first blowing up some codimension $3$ diagonals
corresponding to subsets with more than two elements
to get a space, denoted by $M[n]^{+}$ in [FM], and then
blowing up the proper transform of $V_i$ to get $M[n+1]$,
we only blow up $M[n]\times M$ successively along
the graphs $V_i$. More precisely, 
we first blow up the first graph $V_1$ and 
call the new space ${\cal X}_1$.
Next we blow up
 the proper transformation of the second graph $V_2$ and 
call the new space ${\cal X}_2$. After $n$ steps, we get 
 the space $\tilde M[n]={\cal X}_n$.

 As we are always
 blowing up smooth manifolds along smooth centers, the exceptional divisors
 will be some projective $P^1$ bundles. On the other hand,
 as the graphs intersect each other, one blow-up will definitely affect the
 topological
 types of the other divisors we have already blown up. Thus we do not 
expect the final $n$ smooth exceptional divisors to be simply
 projective $P^1$ bundles. It rather looks like a blown 
up manifold which comes from the $P^1$ bundle by blowing
up a certain low dimensional locus.  Nevertheless, we still have the 
following important proposition.

\proclaim Proposition 4.14. The smooth almost complex manifold $\tilde M[n]$ 
forms a smooth fibre bundle over 
$M[n]$ whose fibres are almost complex four-manifolds all diffeomorphic
to the manifold $M\#n {\overline{\bf CP}}^2$.\par

\noindent{\it Proof}.
 We need to show the natural map $\tilde M[n]\longrightarrow M[n]$
has surjective differential. Away from the exceptional divisors,
this is clear. On the other hand, the exceptional divisors map
to the graphs $V_i$ with surjective differentials, and since
$V_i$ are graphs the differentials of $V_i\longrightarrow M[n]$ are
isomorphisms. Thus 
$\tilde M[n]\longrightarrow M[n]$ is a fibre bundle. The fibres are
diffeomorphic to $M\#{\overline {\bf CP}}^2$ because
the fibres over $M^n(\Delta)$ are just $M$ blown up at $n$
distinct points.

The space
$\tilde M[n]$ does not admit an $S_n$ action;
  the $S_n$ symmetry was broken when we mark the graphs and
determine the order of blowing up using these markings.  
Also notice that the fibres 
actually inherit the almost complex structures from the blow-up process.
 Thus our fibre bundle $\tilde M[n]$
 is a universal fibre bundle parametrizing all the inherited
 almost complex structures
on the   
 manifold $M$ with $n$ points blown up. 

If $M$ is chosen to be a symplectic four-manifold,
 then we can choose the almost complex structure on $M$ to be compatible with
 the given symplectic form. If so, the manifold $M^n$
 becomes a symplectic manifold, too. By performing the symplectic blowing up,
  $\tilde M[n]$ also carries a (non-canonical) symplectic
 structure such that the fibres of the fibration $\tilde M[n]\longrightarrow M[n]$
 are symplectic.  
If one considers the hyperwinding
 family of symplectic forms, one can also construct 
a family of almost complex structures compatible with the given family of
 symplectic forms such that the canonically
 defined sections in these fibre bundles ${\cal X}_n$ are simultaneously
 almost complex(and symplectic)
 with respect to the family of almost complex structures.
 Thus a family version of the Guillemin-Sternberg-McDuff construction gives
 us a family of symplectic structures on the total space $\tilde M[n]$
such that the fibres are symplectic with respect to a family of symplectic
 forms closely related to the original hyperwinding family.

 From now on, we will choose $\tilde M[n]\times S^{b^+-1}$ as the total space 
${\cal X}$ and $M[n]\times S^{b^+-1}$ as the base $B$
and discuss the
 Seiberg-Witten theory on the bundle ${\cal X}\longrightarrow M[n]\times S^{b^+-1}$. 
 As the fibre bundles are not of critical dimensions, 
as mentioned in \S2, to define the Seiberg-Witten invariants we
 need to perturb the Seiberg-Witten equations by the large
self-dual two-forms to avoid reducible solutions. 

In the derivation of the
 wall crossing formula, we do not require the manifold $M$ to be symplectic.
 However, the later application will be for this case. Thus we are free to
 use the symplectic forms to perturb the Seiberg-Witten equations.
 If $M$ is symplectic with $b^+=1$, there are two ways of perturbing the
equations. Given a symplectic two-form, we can perturb the equations either
 by $r\omega$ or $-r\omega$. Letting $r\longrightarrow \infty$, there are two
 chambers and two different invariants. Therefore one should study
the family wall crossing formula to understand how the invariants 
jump. When $b^+$ is greater than one, we can consider 
different winding chambers of self-dual two-forms.
 The family wall crossing formula also makes sense in these cases.

As the fibres are blown up from manifolds diffeomorphic to $M$,  naturally we are 
interested in 
the Spin$^c$ structures which have the following schematic form,
$${\cal L}={\cal L}_0-\sum_i(2m_i+1)E_i$$
 where ${\cal L}_0$ is a 
Spin$^c$ structure on $M$,
$E_i$ is the line bundle which is Poincar\'e dual to 
 the homology class representing the $i-$th exceptional class
 and each 
 $m_i$ is an integer.
Notice that these line bundles $E_i$ give the exceptional line 
bundles of each fibre when they are restricted to each individual fibre.
 The
geometric meaning of these numbers $m_i$ will be clear in the later discussion.   
As before, we still require the dimension of its moduli space to be
non-negative, i.e,
 $$\eqalign {&1/4({\cal L}^2_0-\sum_i(2m_i+1)^2-
(2\tilde\chi+3\tilde \sigma))+4n+(b^+-1)\cr
=
&1/4({\cal L}_0^2-\sum_i(2m_i+1)^2-(2\chi+3\sigma)+n)+4n+(b^+-1)
\cr\geq &0.\cr}$$

Let us analyze the chamber structure of this family.
By Proposition 3.5,
 the period bundle is homotopic to an $S^{b^+-1}$ bundle over 
$M[n]\times S^{b^+-1}$.  Thus by Proposition 3.6, the
chambers are classified by the set of 
fibrewise homotopy classes of the sections
$[M[n]\times S^{b^+-1},{\cal P}]_f$. Usually this set is complicated. We are 
only interested in a subset consisting of those homotopy classes of sections 
whose restriction from
 $M[n]\times S^{b^+-1}$ to $M[n]\times \{pt\}$ is nullhomotopic.
  This subset can be identified with
${\bf Z}_2$ if $b^+=1$
and ${\bf Z}$ if $b^+>1$.
 As the group $\bf Z$ is an infinite cyclic group generated by $1$, we can
 calculate the wall crossing number from the zero winding chamber to the
first winding chamber. Notice that in this case the obstruction class $Obs$
 is the pull back $M[n]\times S^{b^+-1}\longrightarrow S^{b^+-1}$ of the generator 
of $H^{b^+-1}(S^{b^+-1}; {\bf Z})$.  By inserting this special
 obstruction class into the wall
crossing formula we conclude that the calculation can be performed on a certain
projective space bundle over the torus fibration with $M[n]$ as 
its base.
 Notice that it is not necessarily true that the singular locus inside
 $M[n]\times S^{b^+-1}$ is diffeomorphic to $M[n]$. However, 
we will conveniently 
adopt this picture, as the wall 
crossing formula is purely homological.

 In this special case, the wall crossing number calculation is again reduced 
to the calculation of the top power of the hyperplane class over the torus
fibration (when $b_1(M)=0$, the torus fibration reduces to $M[n]$ itself).
As argued in the critical case, this calculation is exactly equivalent to
the calculation of the Segre classes of the Index bundle. 
It is a very complicated calculation if we want to directly apply the family
index theorem. Though there is a slight modification which simplifies the
calculation to certain degree, we are not able to calculate the
final result explicitly. What we can do is to reduce the calculation to
 a problem in algebraic geometry. In the process, it would be manifest 
how the family version of Seiberg-Witten invariants corresponds to the
 similar story on the Gromov side---in terms of algebraic geometry.

We change the strategy as follows. 
Instead of making $\hbox{Coker}_{M[n]}(D)$ a trivial bundle by 
adding its $K$ theory inverse, we do it to $\hbox{Ker}_{M[n]}(D)$. Thus we can 
assume that there is some representative of the Index bundle such that $\hbox{Ker}_{M[n]}(D)$
is trivial, while $\hbox{Coker}_{M[n]}(D)$ absorbs all the topology of the bundles.
Since $\hbox{Coker}_{M[n]}(D)$ forms the inverse bundle of $\hbox{IND}_{M[n]}(D)$ in  the $K$ 
group, the Segre classes of $\hbox{IND}_{M[n]}(D)$ are 
exactly the Chern classes of $\hbox{Coker}_{M[n]}(D)$. In other words, if we want to
calculate the top power of the hyperplane bundle, it is equivalent to 
calculate the top dimensional (in term of the base manifold) Chern class of 
the Cokernel bundle. We would like to demonstrate
 that there is some special candidate for  $\hbox{Coker}_{M[n]}(D)$ which allows us
 to interpret the problem in terms of algebraic geometry.

Before going into the details, let us explain the relationship between the wall 
crossing numbers and the family version Seiberg-Witten
invariants.  In general, to apply the wall crossing formula to calculate the
Seiberg-Witten invariants, we need the information of one specific chamber.
 This is exactly the case one can achieve for the ruled and rational 
surfaces ($b^+=1,B=pt$) and hyperwinding families of the 
K3 surfaces and $T^4$ where we have
some sort of vanishing theorem for the invariants
 due to the existence of
 special metrics of nonnegative scalar curvature.  Then the 
wall crossing number is exactly the Seiberg-Witten
 invariants.  In our case the same
conclusion still holds, but we offer a different argument.
 Instead, we would like to prove that the pure Seiberg-Witten
invariants vanish in the unwinding chamber.

\proclaim Proposition 4.15. (Vanishing result for Seiberg-Witten Invariants).
 The pure family Seiberg-Witten invariants of the $b^+>1$ family
 $\tilde M[n]$ vanish in the zero winding chamber.\par 

\noindent{\it Proof}: If $b^+$ is even, the pure invariants of these 
families are automatically
 zero for dimension reasons. Thus we only need to study the case when $b^+$
is odd. 
Denote the pure Seiberg-Witten invariant in the $i-$th winding chamber
by  $SW({\cal L}, H, i)$; then
 $$\eqalign{WCN( {\cal L},H)=&SW({\cal L},H ,1)-SW({\cal L},H ,0)\cr
 =&SW({\cal L},H ,0)-SW({\cal L},H ,-1.)\cr}$$
   On the other hand, we can use
an orientation reversing diffeomorphism of $S^{b^+-1}$ to map one chamber into
another.  We conclude 
that $$SW({\cal L},H ,-1)=-SW({\cal L},H ,1)$$ 
as changing the base orientation changes the orientation
of $det^+$.
 Therefore 
 $$SW({\cal L},H ,0)=0\qquad\hbox{ and}\qquad
 SW({\cal L},H ,1)=WCN({\cal L},H).$$ 

 We also learn from this proposition that the
wall crossing number calculation gives us the invariants themselves instead
 of their differences. Notice that
 this discussion does not apply to other types of invariants, as other
 types of invariants may have different parity under the orientation reversal 
of the base. On the other hand, the argument offered here is very general.

 Now we would like to 
discuss how to transform our question to a question which 
is closely related to algebraic geometry. For simplicity, let us focus
 on the $b_1=0$ case first. 
 
 Let us look at the family index expression that we want to evaluate;
 it is given by 
$$\int_{\tilde M[n]/M[n]}A_{\tilde M[n]/M[n]}\cdot ch{1\over 2}({\cal L}).
\eqno (4.7)$$
 Even though we evaluate this as a topological number, we will make
use of the
fact that $M$ is an almost complex manifold.
 Writing $\cal L$ as $K^{-1}\otimes C^2$, then (4.7) can 
be recast into 
$$\int_{\tilde M[n]/M[n]}Todd_{\tilde M[n]/M[n]}\cdot
ch(C).\eqno (4.8)$$
 Then the Family index theorem becomes the Grothendieck RR theorem in the
almost complex category.  $C$ can be written as $C_0-\sum_i m_i E_i$
 where $C_0$ is a line bundle on  $M$ and $E_i$ are  line bundles
 dual to the exceptional divisors. 
 
 To evaluate (4.8), 
 we would like to use the Koszul resolution of certain ``sheaves''
in the almost complex category. For this purpose, let us review 
 the usual Koszul resolution in complex geometry. 
 Let $A\subset B$ be a pair
 of complex manifolds such that $A$ is embedded as a complex codimension $d$
 submanifold and is the zero set of a section of a locally free
sheaf $Q$ of rank $d$. Then the structure sheaf ${\cal O}_A$ of $A$, extended by zero,  is a 
skyscraper sheaf on $B$. There is a canonical resolution of ${\cal 
O}_A$
 in terms of locally free sheaves,
$$0\longrightarrow \Lambda^d Q^{\ast}\longrightarrow\cdots \longrightarrow \Lambda^2Q^{\ast}\longrightarrow Q^{\ast}
\longrightarrow{\cal O}_B
\longrightarrow{\cal O}_A\longrightarrow 0,$$
which is usually 
called the Koszul resolution.
 As was remarked at the end of Atiyah-Hirzebruch's paper [AH], this 
resolution goes over 
to the differentiable category as well. Thus it exists even though the manifold is 
almost
complex. Applying to our special case ($d$ is equal to one) we get

$$0\longrightarrow{\cal O}(-E_i)\longrightarrow{\cal O}\longrightarrow{\cal O}_{E_i}\longrightarrow0.$$

 Tensoring this sequence with a locally free sheaf ${\cal N}$
 we get 
$$0\longrightarrow{\cal N}(-E_i)\longrightarrow{\cal N}\longrightarrow{\cal N}_{E_i}\longrightarrow0.$$

 Notice that here we are using the $C^{\infty}$ version of the sheaf
 theory which was sketched in [AH].

Let us first assume $m_i\geq 0$.
In fact, we apply the sequences to different $E_i$ in the reverse 
order. Namely we apply it first to ${\cal N}=C_0-\sum_{i\leq n-1}m_i E_i
-(m_n-1)E_n$. Then we apply it to ${\cal N}=C_0-\sum_{i\leq n-1}m_i E_i
-(m_n-2)E_n$ until the $E_n$ are exhausted. Then we do the similar thing
 to $E_{n-1}$, etc., until all the $E_i$ have been eliminated. In the meantime,
 a sequence of sheaves supported at those $E_i$ will be generated.
 What matters here is that we have found a new representative of ${\cal O}(C_0
-\sum_i m_i E_i)$ in the following form:
$${\cal O}(C_0-\sum_i m_iE_i)\equiv
{\cal O}(C_0)-\sum_{i,j_i;0\leq j_i\leq m_i-1}{\cal O}_{E_i}(
C_0-\sum_{s\leq i-1}m_sE_s-j_iE_i).$$
 As the Chern 
character is an additive homomorphism from the $K-$group to the cohomology
 ring, the family index calculation
 is reduced to the following form.
$$\int_{\tilde M[n]/M[n]} Todd_{\tilde M[n]/M[n]} (ch(C_0)
-\sum_{i,j_i;0\leq j_i\leq m_i-1}{\cal O}_{E_i}
(C_0-\sum_{s\leq i-1}m_s E_s-j_iE_i)).\eqno (4.9)$$

 The first term of (4.9) is interpreted as the index bundle of the $\bar \partial$
operator on the fibres $M\# n\overline{\bf CP}^2$ 
with the coefficient line bundle $C_0$. As $C_0$ is 
pulled back to $\tilde M[n]$ by its projection to $M\times M[n]$, the 
index bundle is a constant bundle on $M[n]$ whose rank is equal to 
 $ind(\bar\partial,C_0)$ on $M$. On the other hand, all the contributions 
of the Chern characters come from the terms which appear with a negative 
sign. As the skyscraper sheaves appear in the Chern character, 
 we can restrict our total space from $\tilde M[n]$ to $E_i$. By 
 the Grothendieck-Riemann-Roch theorem in the differentiable context (see 
[AH] for example), the remaining terms can be cast into the following
 form:
$$-\oplus_{i,j_i;0\leq j_i\leq m_i-1}Ind(\bar\partial,({\cal O}(j_i)))
\otimes (C_0-\sum_{s\leq i-1}m_sE_s)\eqno (4.10)$$
 where the index is taken over the fibres of $E_i$ over $M[n]$ which are
 unions of normal crossing rational curves.

This formula (4.10) can  again be cast into another form. 
To do so requires us 
 to introduce some notation. Inductively, let us consider the
manifold $M\times M[n]$ with the first $i-1$ graphs blown up,
denoted as ${\cal X}_{i-1}$. 
Then the proper transformations of the $i-$th graph 
still forms a section of the projection map from the
blown up manifold ${\cal X}_{i-1}$ 
to $M[n]$. Thus we can identify the restriction of the
vertical tangent bundle (complex two-dimensional) to this section as a 
complex rank two
vector bundle on $M[n]$, which will be denoted by $V_i$. 
Written in terms of $V_i$, the formula (4.10) can be expressed as
$$-\oplus_{i,j_i;0\leq j_i\leq m_i-1}S^{j_i}(V_i^{\ast})\otimes (C_0-
\sum_{s\leq i-1}m_s E_s).\eqno (4.11)$$
where $S^{j_i}(V_i^{\ast})$ denotes the $j_i-$th symmetric power of 
$V_i^{\ast}$.
 Thus the wall crossing number should be expressed as $c_{2n}$ of this
vector bundle. 

If $m_i<0$, then the above discussion is not completely valid.
In this case the Spin$^c$ structure can be rewritten as ${\cal L}_0
+\sum_i(2m_i+1)E_i$ with $m_i$ positive. In this case we still use
 the same exact sequence $$0\longrightarrow {\cal O}(-E_i)\longrightarrow {\cal O}\longrightarrow
 {\cal O}_{E_i}\longrightarrow 0.$$
 But we replace the middle term by the first and the third term.
 Thus the term 
$$\sum_{i,j_i,0\leq j_i\leq m_i}
{\cal O}_{E_i}(C_0+\sum_{s\leq i-1}(m_s+1)E_s+j_iE_i)$$
 still keeps the 
 form similar to the original one. The major difference is that it
appears with a positive sign instead.

 Also $m_i$ is replaced by $m_i+1$.
 The key term is written as $$\oplus_{i,j_i,0\leq j_i\leq m_i}
Ind(\bar\partial,({\cal O}(-j_i)))\otimes (C_0+\sum_{s\leq i-1}(m_s+1)E_s).$$

 We can go through the previous 
argument. This time the line bundles restrict to negative degree line
bundles on the exceptional ${\bf CP}^1$. Recall that the sheaf cohomology
of a negative degree line
bundle on ${\bf CP}^1$ is related to the sheaf cohomology of line bundles of  positive degree
 by the well known Serre duality on curves. We can formally write
$$\eqalign {&[H^0({\bf CP}^1,{\cal O}(-k))]-[H^1({\bf CP}^1,{\cal O}(-k))]\cr
 =&-[H^1({\bf CP}^1,{\cal 
O}(k-2))]^{\ast}
+[H^0({\bf CP}^1,{\cal O}(k-2))]^{\ast}.\cr}$$  
 Applying this formalism to the family of bundles, we still get a series
 of symmetric powers.  The highest power one can get is $m_i+1-2=m_i-1$, the
same index as before. Moreover the Serre duality gives us an extra minus
sign. Thus the answer appears with a negative sign, again compatible with 
the previous answer. Thus we conclude the surprising answer that 
the line bundles ${\cal L}_0+\sum_i(2m_i+1)E_i$, ${\cal L}_0-\sum_i
(2m_i+1)E_i$ share the same family wall crossing number for the section
 invariants.  

At this moment, let us restrict to the case where $M$ is a
 complex manifold (including the case when $ M$ is the K3
surface  or $T^4$). Then 
(4.10) can be interpreted as
 $$-\oplus_{i,j_i;0\leq j_i\leq m_i-1}H^0({E_i}_b,(C_0-\sum_{s\leq i-1}m_sE_s)
\otimes {\cal O}_{{E_i}_b}(j_i)).\eqno (4.12)$$

 It is this expression which has an interesting algebro-geometric meaning.
The first term of the index bundle is nothing else but
 $H^0(M,C_0)$, the space of holomorphic sections of the line bundle $C_0$. 
 The term which is subtracted actually records the $m_i-$jet information
 of the line bundle $C_0$ at the $n$ points of $M$
 parametrized by the point in the
base $M[n]$. There is a natural map from the tautological line bundle over
 $P(H^0(M,C_0))$ to the jet bundle such that this morphism vanishes at a 
point in $P(H^0(M,C_0))\times M[n]$ exactly when the corresponding section 
in $H^0(M,C_0)$ vanishes at these $n$ points with prescribed orders. In other
 words, the counting of such points on the manifold $P(H^0(M,C_0))\times
M[n]$ is equivalent to counting of holomorphic curves (given by the linear
 system corresponding to $C_0$) which have certain singularities along the
$n$
points.   Thus we have seen in this type of examples that the 
algebro-geometric problem is actually encoded in a differential topological invariant.
This gives us very nontrivial evidence that the family version of 
Seiberg-Witten theory should have a certain Gromov-type interpretation.  

Let us come back to the almost complex category. 
In the previous discussion, we have made the assumption $b_1=0$ to
simplify the discussion. Now let us consider the general 
 $b_1\not=0$ case. 

When $b_1=0$, the torus $T^{b_1}$ collapses to
a point and it does not play any role in the formula. 
When $b_1\ne 0$, the key observation is that
 the torus fibration is reduced to $T^{b_1}\times \tilde M[n]$ in this
 special situation. Even though the fibre bundle $\tilde M[n]\longrightarrow
M[n]$ is not
a trivial product bundle, the nontriviality is due to the blow-up process,
which has no
 effect on the $H^1$, nor on
 the flat connections on the manifold $M$. Therefore
 the torus fibration should be the same as the one for $M\times M[n]\longrightarrow
 M[n]$, the trivial product. 

The factorization of the torus fibration
 leads to certain simplifications of the calculations. Previously, the
 Index virtual bundle was split into two parts. The positive part is
 a trivial bundle over the base and the negative part carries all the
 topological information. In the $b_1>0$ case the same conclusion is not valid.
 As the twisting factor $ch(\Omega)$ plays a role,
 the positive factor
 is replaced by a virtual bundle of constant rank. This virtual bundle
 can be written as the formal difference of two terms. The negative part
 is still a trivial bundle. However the positive part is not. If we take
 a closer look at the bundle involved, it is exactly the bundle which 
appears in the family wall crossing formula of ${\cal L}_0$ in the fibre
 bundle $M\times pt$ ($b^+=1$) or $M\times S^{b^+-1}$ ($b^+>1$).
 This vector bundle has the special feature that it is trivial along the
factor $M[n]$.
 On the other hand the negative part of the index bundle is still given by the
 original expression discussed above. It also has the special feature that
 it is trivial along the $T^{b_1}$ factor.
 Thus we get a nice picture out of
 the complicated situation. Namely we are working on a product space
 $T^{b_1}\times M[n]$ and the index bundle is decomposed into terms 
which are pulled back from $T^{b_1}$ and $M[n]$ respectively.

Let us denote the projection maps from $T^{b_1}\times M[n]$ to $T^{b_1}$ and
 $M[n]$ by $\pi_T$ and $\pi_M$ respectively. Then the index bundle
 can be written as 
 $$\hbox{IND}_{M[n]}(D)=\pi_T^{\ast}W_1-\pi_M^{\ast}W_2,\eqno (4.13)$$  where
 $W_1\in K(T^{b_1})$ and $W_2\in K(M[n])$. 

 To calculate the Segre classes of the index bundle 
$S(\hbox{IND}_{M[n]}(D))$, we would like to 
 make use of the product structure described here.
 Instead of doing it directly,
 we can rewrite (4.13) as 
 $$\hbox{IND}_{M[n]}(D)\oplus \pi^{\ast}_MW_2=\pi_T^{\ast}W_1.$$
 As Segre classes are the formal inverses of the Chern classes, they also
 satisfy the Whitney formula. Thus we have $S(\hbox{IND}_{M[n]}(D))\cdot S(\pi_M^{\ast}W_2)
=S(\pi_T^{\ast}W_1)$. Multiplying both sides by 
$C(\pi_M^{\ast}W_2)$ and making use of the $S\cdot C=1$ property we
 finally get $$S(\hbox{IND}_{M[n]}(D))=S(\pi^{\ast}_TW_1)\cdot C(\pi_M^{\ast}W_2).$$

 If we are interested in the general terms of the left hand side,
 it is expressed in terms of the combinations of the Chern classes and
 Segre classes of $W_1$ and $W_2$ in a cumbersome way. However the result we are
interested in is the top Segre class (where the word ``top'' means 
 the top dimension of
 $T^{b_1}\times M[n]$). In this case the formula simplifies and it is
 given by the top Segre class of $W_1$ and the top Chern class of $W_2$. 
 We also notice that the top Segre class of $W_1$ is nothing else than
 the wall crossing number of ${\cal L}_0$ in the critical product family.
 Thus, we have obtained an amazing result relating the family wall crossing
 formula of the fibre bundle $\tilde M[n]$(or $\tilde M[n]\times S^{b^+-1}$)
 and the family wall crossing formula of $M$ or $M\times S^{b^+-1}$.
The following is the main theorem of this subsection.

\proclaim Theorem 4.16. Let $\delta_{M[n]}({\cal L})$ be the wall crossing number
 of the consecutive chambers of the family $\tilde M[n]$ or $\tilde M[n]\times
 S^{b^+-1}$.
  Let ${\cal L}$ be written as ${\cal L}_0+\sum_i(2m_i+1)E_i$
 with $m_i\geq 0$; then $\delta_{M[n]}({\cal L})$ can be written as the
 following expression in terms of $\delta({\cal L}_0)$, the wall crossing
 number of the consecutive chambers of the family $M\times pt$ or 
$M\times S^{b^+-1}$,
 $$\delta_{M[n]}({\cal L})=\int_{M[n]}c_{2n}(\oplus_{i,j_i;0\leq j_i\leq m_i-1}
S^{j_i}(V_i^{\ast})\otimes (C_0-\sum_{s\leq i-1}m_sE_s))\cdot
 \delta({\cal L}_0), \eqno (4.14)$$ 
 where ${\cal L}_0=2C_0-K_M$. 
  If $2m_i+1$ is negative in the formula of ${\cal L}$, then its wall crossing
 number is the same as the one with $-(2m_i+1)$ replacing $2m_i+1$. Namely,
 there is an explicit ${\bf Z}_2$ symmetry among the Spin$^c$ structures which
 preserves the wall crossing number.\par

When $b_1(M)=0$, the wall crossing number of the critical family becomes
 $\pm 1$ and the formula goes back to the expression studied earlier.
 We want to emphasize that only in the $b_1\not=0$ case does the dependence
 on $\delta({\cal L}_0)$ become explicit.

 The theorem is not only valid for the pure family Seiberg-Witten invariant,
 it is valid for the invariant involving $H_1(M;{\bf Z})$ as well,
 as long as $\delta({\cal L}_0)$ is understood to be the wall crossing
 number of the corresponding invariant of the critical family.

 On the other hand, if we are interested in the wall crossing numbers of the
 mixed invariants, it is easy to derive the corresponding wall crossing 
formula. Let $\eta$ be the cohomology class in $H^{\ast}(M[n];{\bf Z})$ whose
 insertion defines the mixed Seiberg-Witten invariants. Then the wall 
crossing formula is of almost the same form except that 
 $\int_{M[n]}c_{2n}(\cdots)$ is replaced 
by $\int_{M[n]}\eta\cdot c_{2n-deg\eta}
(\cdots)$. In particular,  the wall crossing formula for the mixed invariant
 with the fundamental class $[M[n]]\in H^{top}(M[n];{\bf Z})$ inserted is the
 same as the wall crossing formula $\delta({\cal L}_0)$.
 
 There are several things we learn from the derivation
 of the wall crossing formula for these special families
constructed out of the Fulton-MacPherson spaces. 
 First, the lengthy formula in terms of all the Chern characters may not
 be so useful in deriving the formula. One may need to figure out 
certain tricks (in our case, the Koszul resolution of sheaves in the smooth
 category) to simplify the calculation. As the final answer depends on the
 topology of the fibre bundle, one should not expect the answer to be 
 as universal  as in the critical case. To calculate the formula
 for  other interesting families is an interesting question on its own.
 
 Second, as the wall crossing formula satisfies the amazing factorization
 property, one may suspect that the Seiberg-Witten invariants
 themselves also satisfy similar constraints. Namely, they are related
 to the family invariants of the critical families in an explicit way.
 In fact, the family blow-up formula [Liu] supports this. Using the
 family blow-up formula, one is able to check that 
the family Seiberg-Witten 
invariants satisfy the constraint predicted here. For the
 manifolds $M$ with $b^+-1$  positive
 and even, the result could be obtained by the wall crossing formula and
 the vanishing result stated earlier. However, the family blow-up 
formula gives us an independent derivation valid for all cases (including
 the $b^+=1$ case). The compatibility of the blow-up formula and the
 wall crossing formula gives us a commutative diagram between blowing up
 and crossing chambers. Let us begin from the fibre bundle $M\times 
M[n]\times S^{b^+-1}$. One can either move from the zero winding chamber
 to the nonzero winding chamber and then blow up to $\tilde M[n]\times 
S^{b^+-1}$. Or one can first blow up $n$ times
 to $\tilde M[n]\times S^{b^+-1}$
 and then go across the walls to the nonzero winding chambers.
 The family blow-up formula also explains why the vanishing result
 holds for the zero winding chamber. Once the pure invariants of the
 zero winding chamber of the critical families are known to be zero, the
 others follow as well by the family blow-up formula. For the details
  consult [Liu]. \par

 Finally, let us give the explicit formulas when $n=1$. 
 When $n=1$, $M[n]$ is simply $M$, and $V_1$ is just the tangent bundle.
 Consider the case $m_1=2$; then by (4.14)
 $$\eqalign{\delta_M({\cal L})=&\delta({\cal  L}_0)\int_{M} c_2(
 (S^0(TM)\oplus S^1(TM) )\otimes C_0)\cr
 =&\delta({\cal  L}_0 )(3C_0\cdot C_0 +2C_0\cdot c_1(M) + c_2(M))\cr}$$
 
 When $M$ is the four-torus, both $c_1(M)$ and $c_2(M)$ are trivial,
 so by Corollary 4.11,
 $$\eqalign{\delta_M({\cal L})=&\delta({\cal  L}_0 )3C_0\cdot C_0\cr
 =&{1\over 2}(C_0\cdot C_0+1) 3C_0\cdot C_0.\cr} \eqno (4.15)$$
 
 In [BL2], for a hyperk\"ahler family of four-tori,
  the number of $1-$nodal curves in a primitive class $C_0$, 
 is exactly given by (4.15). This is again evidence supporting
`SW=TGW'.

\bigskip
\bigskip 
\noindent{\bf References}.
\medskip
\item{} [AH] M. Atiyah and F. Hirzebruch, The Riemann-Roch theorem for 
analytic embeddings, Topology 1 (1962) 151-166.
\item{} [BL1] J. Bryan and C. Leung, The enumerative geometry of K3 surfaces
and modular forms, J. Amer. Math. Soc. 13 (2000), No.2, 371-410. 
\item{} [BL2] J. Bryan and C. Leung, Generating functions for the 
number of curves on Abelian surfaces, Duke Math. J. 99 (1999), No.2, 311-328.
\item{} [D1] S. Donaldson, The Yang-Mills invariants of four-manifolds, 
{\it Geometry of low-dimensional manifolds, I (Durham, 1989)}, 5-40, London Math. Soc. 
Lecture Note Ser., 150, {\it Cambridge Univ. Press, Cambridge, 1990}
\item{} [D2] S. Donaldson,  The Seiberg-Witten equations
and $4$-manifold topology,  Bull. AMS 33 (1996), no.1, 45-70.
\item{}[DK] 
S. Donaldson and P. Kronheimer,
{\it The geometry of four-manifolds}. Oxford Mathematical Monographs. Oxford
Science Publications. {\it The Clarendon Press, Oxford University Press, New
York, 1990}.
\item{}[F] 
W. Fulton,
{\it Intersection theory}. Ergebnisse der Mathematik und ihrer Grenzgebiete (3)
[Results in Mathematics and Related Areas (3)], 2. {\it Springer-Verlag,
Berlin, 1984}.                     
\item {}[FS] R. Fintushel and R. Stern, Knots, links, and 4-manifolds, Invent. Math.
134 (1998), No.2, 363-400.
\item {}[FM] W. Fulton and R. MacPherson, A compactification
of configuration spaces, Ann. of Math. 139 (1994)183-225.
\item {}[FMo] R. Friedman and J. Morgan, Algebraic
surfaces and Seiberg-Witten invariants, J. Algebraic Geom. 6 (1997), No.3, 445-479.
\item{} [Ge] H. Geiges, symplectic couples on 4-manfolds, Duke. Math. Journal
85 (1996), 701-711.
\item{} [Gr] M. Gromov, Pseudo-holomorphic curves in symplectic 
manifolds, Invent. Math. 82 (1985), 307-347.
\item{} [GS] V. Guillemin and S. Sternberg, Birational equivalence in the symplectic category, Invent. Math. 97 (1989), No.3, 485-522. 
\item{} [IP] E-N. Ionel and T. Parker, The Gromov invariants of Ruan-Tian 
and Taubes, Math. Res. Lett. 4 (1997), No.4, 521-532.
\item{} [K] P. Kronheimer, Some nontrivial families of symplectic structures,
preprint.
\item{} [KM] P. Kronheimer and T. Mrowka, The genus
of embedded surfaces in the projective plane,
Math. Res. Lett. 1 (1994), 797-808.
\item{} [Liu] A. Liu, Family blow up formula and nodal curves
in K\"ahler surfaces, preprint.
\item{} [LL1] T. J. Li and A. Liu, General wall crossing formula,
Math. Res. Lett. 2 (1995), No. 6, 797-810.
\item{} [LL2] T. J. Li and A. Liu, On the equivalence between
SW and GT in the case $b^+=1$, Internat.
Math. Res. Notices 1999, no. 7, 335--345.
\item{} [LL3] T. J. Li and A. Liu, Symplectic four-manifolds with torsion
canonical classes, in preparation. 
\item{} [M] D. McDuff, Lectures on
Gromov invariants for symplectic 4-manifolds, in {\it Gauge theory and Symplectic Geometry
(Montreal, 1995)}, 175-210, Kluwer Acad. Publ. 1997.
\item{} [M2]. D. McDuff, The structure of rational
and ruled symplectic $4-$manifold, J. Amer. Math. Soc. v.1. no.3, (1990), 679-710.
\item{}[R] 
Y. Ruan,
Virtual neighborhoods and pseudo-holomorphic curves. Proceedings of 6th
Gökova Geometry-Topology Conference. Turkish J. Math. 23 (1999), no. 1,
161--231.
\item{} [RT] Y. Ruan and G. Tian, A mathematical
theory of quantum cohomology,  J. Diff. Geom. 42 (1995), no.2, 259-367.
\item{}[Rub] D. Ruberman, An obstruction to smooth isotopy in dimension 4,
Math. Res. Lett. 5 (1998), No.6, 743-758.
\item{} [S] N. Steenrod, {\it Topology of fibre bundles ,
Princeton University Press, Princeton, 1951}.
\item{}[SW] N. Seiberg and E. Witten,
Electric-magnetic duality, monopole
condensation, and confinement in $N=2$ supersymmetric Yang-Mills theory,
Nuclear Physics B, 426 (1994), 19-52.
\item{} [T1] C. H. Taubes, The Seiberg-Witten
invariants and symplectic forms, Math. Res. Letters 1 (1994), 809-822.
\item{} [T2] C. H. Taubes,  SW$\Rightarrow$Gr: 
From the Seiberg-Witten equations to pseudo-holomorphic curves, 
Jour. Amer. Math. Soc. 9(1996), 845-918.
\item{} [T3] C. H. Taubes, Counting pseudo-holomorphic submanifolds
in dimension four, J. Diff. Geom. 44 (1996), no. 4, 818-893.
\item{} [T4] C. H. Taubes, Gr$\Rightarrow$SW: From pseudo-holomorphic curves to
the Seiberg-Witten solutions, J. Diff. Geom. 51 (1999), No.2, 203-334.
\item{} [T5] C. H. Taubes, Gr=SW, Counting curves and connections, J. Diff. Geom. 52 (1999), No.3, 
453-609.
\item{} [YZ] S.T. Yau and E. Zaslow, BPS states, string duality, and nodal 
curves on $K3$, Nuc. Phys. B 471 (1996), no.3, 503-512.
\item{}[W] 
E. Witten, Monopoles and four-manifolds, Math. Res. Lett. 1 (1994),
769-796.

\bigskip

\item{}Department of Mathematics, Yale University, CT 06520 
\item{}tli@math.yale.edu
\item{}{\it Current address}:
\item{}Department of Mathematics, Princeton University, Princeton, NJ 08544
\item{}tli@math.princeton.edu
\medskip
\item{}Department of Mathematics, MIT, Cambridge, MA 02139
\item{}akliu@math.mit.edu
\item{}{\it Current address}:
\item{}Department of Mathematics, UC Berkeley, Berkeley, CA 94720     
\item{}akliu@math.berkeley.edu

\end